\documentclass[11pt]{amsart}
\usepackage{latexsym,graphicx}

\numberwithin{equation}{section}
\theoremstyle{plain}

\theoremstyle{remark}

\theoremstyle{definition}

\newcommand{\D}{{\mathcal D}}
\newcommand{\E}{\mathcal E}

\newcommand{\G}{{\mathcal G}}

\newcommand{\K}{{\mathcal K}}
\renewcommand{\L}{{\mathcal L}}
\newcommand{\M}{{\mathcal M}}
\newcommand{\N}{\mathbb N}

\newcommand{\V}{{\mathcal V}}

\newcommand{\dist}{\operatorname{dist}}

\newcommand{\fp}{\operatorname{FP}}

\newcommand{\id}{\operatorname{id}}

\newcommand{\Int}{\operatorname{Int}}

\renewcommand{\span}{\operatorname{span}}

\newcommand{\supp}{\operatorname{Supp}}

\def\half{{1 \over 2}}
\def\la{\lambda}

\newcommand{\oa}{\overrightarrow}

\newcommand{\ol}{\overline}

\def\XXint#1#2#3{{\setbox0=\hbox{$#1{#2#3}{\int}$}
      \vcenter{\hbox{$#2#3$}}\kern-.5\wd0}}

\address{ Dept. of Mathematics, Rice  University, 6100 Main St., Houston, 77005 TX, U.S.A.
\\ {\sl E-mail address:}  {\bf harvey@rice.edu}}

\address{ Department of Mathematics, Stony Brook University, Stony Brook, NY 11790, U.S.A.
\\ {\sl E-mail address:}  {\bf blaine@math.sunysb.edu}}

\begin{document}

\def\cal{\mathcal}

\font\tpt=cmr10 at 12 pt
\font\fpt=cmr10 at 14 pt

\font \fr = eufm10



\overfullrule=0in

\def\boxit#1{\hbox{\vrule
 \vtop{%
  \vbox{\hrule\kern 2pt %
     \hbox{\kern 2pt #1\kern 2pt}}%
   \kern 2pt \hrule }%
  \vrule}}

  \def\harr#1#2{\ \smash{\mathop{\hbox to .3in{\rightarrowfill}}\limits^{\scriptstyle#1}_{\scriptstyle#2}}\ }

 \def\GG{{{\bf G} \!\!\!\! {\rm l}}\ }

\def\GL{{\rm GL}}

\def\bll{I \!\! L}

\def\bra#1#2{\langle #1, #2\rangle}
\def\bbf{{\bf F}}
\def\bbj{{\bf J}}
\def\Jtn{{\bbj}^2_n}  \def\JtN{{\bbj}^2_N}  \def\JoN{{\bbj}^1_N}
\def\jt{j^2}
\def\jtx{\jt_x}
\def\Jt{J^2}
\def\Jtx{\Jt_x}
\def\bpp{{\bf P}^+}
\def\bpt{{\wt{\bf P}}}
\def\fsh{$F$-subharmonic }
\def\mo{monotonicity }
\def\jet{(r,p,A)}
\def\ss{\subset}
\def\sse{\subseteq}
\def\half{\hbox{${1\over 2}$}}
\def\smfrac#1#2{\hbox{${#1\over #2}$}}
\def\oa#1{\overrightarrow #1}
\def\dim{{\rm dim}}
\def\dist{{\rm dist}}
\def\codim{{\rm codim}}
\def\deg{{\rm deg}}
\def\rank{{\rm rank}}
\def\log{{\rm log}}
\def\Hess{{\rm Hess}}
\def\Hessyp{{\rm Hess}_{\rm SYP}}
\def\trace{{\rm trace}}
\def\tr{{\rm tr}}
\def\max{{\rm max}}
\def\min{{\rm min}}
\def\span{{\rm span\,}}
\def\Hom{{\rm Hom\,}}
\def\det{{\rm det}}
\def\End{{\rm End}}
\def\Sym{{\rm Sym}^2}
\def\diag{{\rm diag}}
\def\pt{{\rm pt}}
\def\Spec{{\rm Spec}}
\def\pr{{\rm pr}}
\def\Id{{\rm Id}}
\def\Grass{{\rm Grass}}
\def\Herm#1{{\rm Herm}_{#1}(V)}
\def\arr{\longrightarrow}
\def\supp{{\rm supp}}
\def\Link{{\rm Link}}
\def\Wind{{\rm Wind}}
\def\Div{{\rm Div}}
\def\vol{{\rm vol}}
\def\foral{\qquad {\rm for\ all\ \ }}
\def\fpsh{{\cal PSH}(X,\f)}
\def\Core{{\rm Core}}
\def\dis{f_M}
\def\Re{{\rm Re}}
\def\rn{\bbr^n}
\def\pp{\cp^+}
\def\plp{\cp_+}
\def\Int{{\rm Int}}
\def\cix{C^{\infty}(X)}
\def\Gr#1{G(#1,\rn)}
\def\Symn{{\Sym(\rn)}}
\def\SymN{{\Sym(\bbr^N)}}
\def\Gpn{G(p,\rn)}
\def\fd{{\rm free-dim}}
\def\SA{{\rm SA}}
 \def\cd{{\cal C}}
 \def\cdt{{\widetilde \cd}}
 \def\cm{{\cal M}}
 \def\cmt{{\widetilde \cm}}

\def\Theorem#1{\medskip\noindent {\bf THEOREM \bf #1.}}
\def\Prop#1{\medskip\noindent {\bf Proposition #1.}}
\def\Cor#1{\medskip\noindent {\bf Corollary #1.}}
\def\Lemma#1{\medskip\noindent {\bf Lemma #1.}}
\def\Remark#1{\medskip\noindent {\bf Remark #1.}}
\def\Note#1{\medskip\noindent {\bf Note #1.}}
\def\Def#1{\medskip\noindent {\bf Definition #1.}}
\def\Claim#1{\medskip\noindent {\bf Claim #1.}}
\def\Conj#1{\medskip\noindent {\bf Conjecture \bf    #1.}}
\def\Ex#1{\medskip\noindent {\bf Example \bf    #1.}}
\def\Qu#1{\medskip\noindent {\bf Question \bf    #1.}}
\def\Exercise#1{\medskip\noindent {\bf Exercise \bf    #1.}}

\def\HoQu#1{ {\AAA T\BBB HE\ \AAA H\BBB ODGE\ \AAA Q\BBB UESTION \bf    #1.}}

\def\pf{\medskip\noindent {\bf Proof.}\ }
\def\qed{\hfill  $\vrule width5pt height5pt depth0pt$}
\def\equdef{\buildrel {\rm def} \over  =}
\def\qedqed{\hfill  $\vrule width5pt height5pt depth0pt$ $\vrule width5pt height5pt depth0pt$}
\def\mathqed{  \vrule width5pt height5pt depth0pt}

\def\V{W}

\def\df{d^{\phi}}
\def\hk{\_{\rm l}\,}
\def\n{\nabla}
\def\w{\wedge}

\def\cu{{\cal U}}   \def\cc{{\cal C}}   \def\cb{{\cal B}}  \def\cz{{\cal Z}}
\def\cv{{\cal V}}   \def\cp{{\cal P}}   \def\ca{{\cal A}}
\def\cw{{\cal W}}   \def\co{{\cal O}}
\def\ce{{\cal E}}   \def\ck{{\cal K}}
\def\ch{{\cal H}}   \def\cm{{\cal M}}
\def\cs{{\cal S}}   \def\cn{{\cal N}}
\def\cd{{\cal D}}
\def\cl{{\cal L}}
\def\cp{{\cal P}}
\def\cf{{\cal F}}
\def\ccr{{\cal  R}}

\def\gerG{{\fr{\hbox{g}}}}
\def\gerB{{\fr{\hbox{B}}}}
\def\gerR{{\fr{\hbox{R}}}}
\def\p#1{{\bf P}^{#1}}
\def\vf{\varphi}

\def\wt{\widetilde}
\def\wh{\widehat}

\def\and{\qquad {\rm and} \qquad}
\def\arr{\longrightarrow}
\def\ol{\overline}
\def\bbr{{\mathbb R}}\def\bbh{{\mathbb H}}\def\bbo{{\mathbb O}}
\def\bbc{{\mathbb C}}
\def\bbq{{\mathbb Q}}
\def\bbz{{\mathbb Z}}
\def\bbp{{\mathbb P}}
\def\bbd{{\mathbb D}}

\def\a{\alpha}
\def\b{\beta}
\def\d{\delta}
\def\e{\epsilon}
\def\f{\phi}
\def\g{\gamma}
\def\k{\kappa}
\def\la{\lambda}
\def\o{\omega}

\def\s{\sigma}
\def\x{\xi}
\def\z{\zeta}

\def\D{\Delta}
\def\L{\Lambda}
\def\G{\Gamma}
\def\O{\Omega}

\def\bd{\partial}
\def\bdf{\partial_{\f}}
\def\lag{Lagrangian}
\def\psh{plurisubharmonic }
\def\ph{pluriharmonic }
\def\pph{partially pluriharmonic }
\def\omp{$\omega$-plurisubharmonic \ }
\def\ffl{$\f$-flat}
\def\PH#1{\widehat {#1}}
\def\lloc{L^1_{\rm loc}}
\def\dbar{\ol{\partial}}
\def\lp{\Lambda_+(\f)}
\def\lpp{\Lambda^+(\f)}
\def\bo{\partial \Omega}
\def\Ob{\overline{\O}}
\def\fc{$\phi$-convex }
\def\PSH{{ \rm PSH}}
\def\SH{{\rm SH}}
\def\totr{ $\phi$-free }
\def\BM{\lambda}
\def\Der{D}
\def\CH{{\cal H}}
\def\RH{\overline{\ch}^\f }
\def\pconv{$p$-convex}
\def\MA{MA}
\def\lagpsh{Lagrangian plurisubharmonic}
\def\hermsk{{\rm Herm}_{\rm skew}}
\def\PSHl{\PSH_{\rm Lag}}
 \def\ppsh{$\pp$-plurisubharmonic}
\def\fp{$\pp$-plurisubharmonic }
\def\fh{$\pp$-pluriharmonic }
\def\Symn{\Sym(\rn)}
 \def\ci{C^{\infty}}
\def\USC{{\rm USC}}
\def\fa{{\rm\ \  for\ all\ }}
\def\ppc{$\pp$-convex}
\def\cpt{\wt{\cp}}
\def\ft{\wt F}
\def\ob{\overline{\O}}
\def\Be{B_\e}
\def\K{{\rm K}}

\def\M{{\bf M}}
\def\N#1{C_{#1}}
\def\ds{Dirichlet set }
\def\dir{Dirichlet }
\def\Fa{{\oa F}}
\def\TR{{\cal T}}
 \def\LAG{{\rm LAG}}
 \def\ISO{{\rm ISO_p}}
 \def\Span{{\rm Span}}
 \def\sk{{\rm skew}}
 \def\sym{{\rm sym}}
 \def\cn{\bbc^n}
 \def\Herm{{\rm Herm}}
 \def\cpp{\cp^+}
 \def\cpm{\cp_+}
 \def\tobd{\bigr|_{\bo}}
 \def\rest#1{\bigr|_{#1}}

\def\AA{1}
\def\BB{2}
\def\CC{3}
\def\DD{4} 
\def\DDD{5}
\def\EE{6}
\def\FF{7}
\def\GGG{8}
\def\HH{9}
\def\II{10}
\def\JJ{11}

\vskip .4in

\def\Bry{1}
\def\CRA{2}
\def\CIL{3}
\def\Har{4}
\def\DDd{5}
\def\PUP{6}
\def\DDR{7}
\def\Survey{8}
\def\AE{9}
\def\LA{10}

\def\PSF{{\mathcal PSH}}

\def\Ext{{\rm Ext}}
\def\E{E}
\def\bL{{\bf \Lambda}}
\def\hn{{\bbh^n}}
\def\Lag{{\rm Lag}}
\def\bu{B^\uparrow}

\font\headfont=cmr10 at 14 pt

\vskip .1in


\title[PLURIHARMONICS IN GENERAL POTENTIAL THEORIES]
{
PLURIHARMONICS IN GENERAL POTENTIAL THEORIES
}

\date{\today}
\author{ F. Reese Harvey \  \and\   H. Blaine Lawson, Jr.
}
 \thanks{The second author was partially supported by the NSF.
}

\maketitle

\centerline{\bf Abstract}
  \font\abstractfont=cmr10 at 10 pt
  
  {{\parindent= .2in\narrower \noindent
  
  The general purpose of this paper is to investigate the notion of ``pluriharmonics''  
  for the general potential theory associated to a convex cone subequation $F\ss \Symn$.
  For such $F$ there exists a maximal linear subspace $E\ss F$, called the {\sl edge},  and $F$ decomposes
  as $F=E \oplus F_0$. The {\sl pluriharmonics} or {\sl edge functions} are $u$'s with $D^2u \in E$.
   Many subequations $F$ have the same edge $E$, but there is a
  unique smallest such subequation.  These are the focus of this investigation.
  Structural results are given.  Many examples are described, and a classification of 
  highly symmetric cases is given.  Finally, the relevance of edge functions to the solutions
  of the Dirichlet problem is established.

}}

\vskip .3in
\centerline{\bf Table of Contents}

\noindent
 \AA.     Introduction.   
 
\noindent
\BB.    Preliminaries -- The Edge and the Span of a Subequation.

\noindent
 \CC.  Extremely Degenerate Versus Completeness.

\noindent
 \DD.    The Supporting Subequation.

\noindent
 \DDD.   Minimal  Subequations.

\noindent
 \EE.    Edge Functions -- Pluriharmonics.

\noindent
 \FF.   Further Discussion of Examples.

\noindent
 \GGG.   Classifying the Minimal Subequations.

\noindent
 \HH.  An Envelope Problem for Minimal Subequations.
 

\vfill\eject


\medskip
\noindent{\headfont \AA.\  Introduction.}

This paper is concerned with the edge $E$ of a  convex cone subequation $F\ss \Symn$,
obtained from the decomposition
$$
F\ =\ E \oplus F_0
\eqno{(\AA.1)}
$$
into a  vector subspace $E\ss\Symn$ and a cone   $F_0 \ss E^\perp$,
called the {\sl reduced constraint set}, which contains no lines.
(See \S  2 for definition of {\sl subequation}.)
The interest in  the edge $E$ is that it gives us a 
notion of {\sl pluriharmonics}, or {\sl edge functions},  for the potential theory
associated to the subequation $F$.
These edge functions $u$ are defined by $D^2u\in E$.

The edge $E$ of a subequation is, in a sense, a crude invariant, since many 
subequations $F$ have the same edge $E$.  However, there is a canonical
choice for $F$ completely determined by $E$, namely $E+\cp$, where $\cp \equiv \{A : A\geq0\}$.
This is a subequation with edge $E$, and it must be contained in all other subequations
with edge $E$ (since by definition we always have $\cp\ss F$).
A large part of this paper is devoted to studying and classifying these subequations
$E+\cp$ for which the edge $E$ is a determinative invariant.  It is for these {\bf minimal
 subequations}  with edge $E$ (Def. \DDD.2) that the pluriharmonics or edge functions play the most important role.

Let's look at some examples.  The simplest is just $\cp$ itself.
Here the edge $E=\{0\}$ and the $\cp$-subharmonics are the convex
functions.  The edge functions  are those  $u$ with $D^2 u\equiv 0$,
that is, the affine functions.

The other orthogonally invariant edge is $E=\{A : \tr A = 0\}$, and we have
$\cp+E =\{A : \tr A\geq 0\} \equiv \D$.  Here the subsolutions are the classical subharmonics,
and the edge functions are just the harmonics.

The examples become  more interesting if we look at U$(n)$-invariant edges in 
$\cn = (\bbr^{2n}, J)$. The reduction into irreducibles is
$$
\Sym(\bbr^{2n}) \ =\ (\bbr\cdot {\rm Id}) \oplus \Herm^\sym_0\oplus \Herm^\sk,
$$
 $\Herm^\sym_0 = \{A : AJ=JA \ {\rm and}\ \tr A=0\}$ and 
 $\Herm^\sk = \{A : AJ=-JA \}$.
Here there are two new edges.  

The first is where we set $E=\Herm^\sk$.  This gives the complex Monge-Amp\`ere subequation:
$$
\cp_\bbc \ =\ E+\cp \ =\ \{A : A_\bbc = A-JAJ \geq 0\}.
$$
The subsolutions are the  plurisubharmonics, and the edge functions (or pluriharmonics) are
the classical pluriharmonic functions in complex analysis.

The other is where  $E=\Herm^\sym_0$. This subequation is rather new.
$$
\cp(\Lag)\ =\ E+\cp\ =\ \left \{A : \tr \left (A\bigr|_W \right)\geq0 \ \ \text{for all Lagrangian planes } W\right \}.
$$
The subsolutions are the {\sl Lagrangian plurisubharmonics} which were studied in [\LA].  In this case
the edge functions are certain  quadratic functions.

If one now looks for Sp$(n)\cdot$Sp(1)-invariant edges, there are many interesting examples.
In fact, in Chapter \FF\  a wide variety of  subequations are given.
The reader might enjoy this section.

One interesting result concerning these  minimal  subequations  is 
that  they can be used to characterize the dual subharmonics (negatives of the superharmonics).
(In this paper, degree-2 means degree $\leq2$.)

 \Theorem{\EE.4}
{\sl
Let $\cp^+=E+\cp$ be a minimal subequation.
The following conditions on an
upper semi-continuous function $u$ are equivalent.
$$
\begin{aligned}   
&\text{(1)\ \ $u$ is dually $\cpp$-subharmonic}. \\ 
&\text{(2)\ \ $u$ is ``sub'' the edge functions}. \\
 &\text{(3)\ \ $u$ is locally ``sub'' the edge functions}. \\
 &\text{(4)\ \ $u$ is locally ``sub'' the degree-2 polynomial edge functions}.
\end{aligned}
$$
}
\noindent
Here the notion of subharmonic is taken in the viscosity sense (see below).
We note that $u$ is ``sub'' a function $v$ if  for each compact set $\ob$, one has
  $u\leq v$ on $\bo  \Rightarrow u\leq v$ on $\ob$.

For example, when $F=\cp$ (the first example above), this theorem says that
the $\cpt$-subharmonic functions (the negatives of supersolutions for the real Monge-Amp\`ere equation)
are characterized by being ``sub'' the affine functions.

When $F=\cp_\bbc$, the $\cpt_\bbc$-subharmonic functions (the negatives of supersolutions for the complex Monge-Amp\`ere equation) are characterized by being ``sub'' the 
standard pluriharmonics in complex analysis.  In fact the degree-2 polynomial pluriharmonics 
will do.

When $F=\cp(\Lag)$, the $\wt \cp(\Lag)$-subharmonic functions (the negatives of  supersolutions for the Lagrangian
 Monge-Amp\`ere equation [\LA]) are characterized by being ``sub''  an explicit family of quadratic polynomials.

Chapter \BB \ of this paper lays down the fundamental notions of the edge $E$, the span  $S$ and the reduced
constraint set $F_0$, which appears in (\AA.1) above.  By definition $S$ is the linear span of $F_0$, but
it could also be defined as the orthogonal complement of $E$. This chapter then looks 
 at  geometrically defined equations, which give many interesting examples.

Chapter \CC\  introduces  a further refinement of the structure of convex cone subequations.
These equations are partitioned into two classes.  The first consists of subequations
which are {\sl extremely degenerate} (see Definition \CC.1).  A basic example is the Laplacian
on $\bbr^k$, considered as a subequation on $\rn$ where $n>k$.
The second class consists of those subequations which are {\sl dimensionally complete} 
 (see Definition \CC.1$^*$). This means that all the variables in $\rn$ are required to define
 the subequation.  These complementary concepts are determined entirely by properties
of the edge, or equivalently the span, of $F$.  For example,
$$
\text{$F$ is complete $\quad\iff\quad E\cap \cp= \{0\}
    \quad\iff\quad S\cap (\Int \cp) \neq \emptyset$.}
$$
There are several other equivalent criteria; see
Propositions \CC.5 and \CC.10.
A pair of orthogonal  subspaces $E$ and $S$ of $\Symn$, which satisfy 
these criteria for $E$ and $S$, is called an {\bf edge-span pair}.

In Chapter \DD\ the structure of the subequation $F$ is further illuminated by proving 
that there exists a unique subspace $W\ss\rn$, called the {\bf support of $F$}, with the 
property that 
$$
F\ =\ \Sym(W)^\perp \oplus F_1 \qquad{\rm where} \ \ F_1\ss\Sym(W)
\eqno{(\AA.2)}
$$
is a complete convex cone subequation in $W$. This $F_1$ is called the {\bf the supporting subequation
of $F$}.

Combining (\AA.1) and (\AA.2) give the decomposition
$$
F \ =\ \Sym(W)^\perp \oplus E_1 \oplus F_0
\eqno{(\AA.3)}
$$
where the edge of $F$ is $E=\Sym(W)^\perp \oplus E_1$ and $E_1$ is the
edge of the supporting subequation $F_1$.

In Chapter \DDD\ our minimal subequations are defined and discussed.
We start  with any  basic edge $E\ss\Symn$,  i.e, one which satisfies
$E\cap\cp=\{0\}$.  Then (Lemma \DDD.1)  the sum 
$$
\cpp \ \equiv \ E+\cp \quad\text{is a subequation, and it has edge $E$}.
$$
(By the edge criteria $\cpp$ is then complete.)
Such subequations are called {\bf minimal} and are the main focus of this paper.
These subequations have many special properties.
Theorem \DDD.4 mentions eight of them, while Theorem \DDD.5 claims that 
any one of these eight properties implies the subequation is minimal.

In Chapter \FF\  many  examples of minimal  subequations are given.

In Chapter \GGG\  we classify all the minimal subequations which are invariant under
the compact group $G$ where $G= {\rm O}_n, {\rm U}_n, {\rm Sp}_n\cdot {\rm Sp}_1,
 {\rm Sp}_n$ and  ${\rm Sp}_n\cdot {\rm S}^1$ (all acting on their fundamental 
 representation spaces).

We note that the general definition of $F$-plurisubharmonics is based on viscosity theory [\CIL], [\CRA]. 
The reader is referred to  [\DDd], [\DDR] or [\Survey] for definitions and properties.

In Chapter \EE\ the generalized  pluriharmonics (or edge functions) are introduced  from a viscosity
point of view, and basic properties are discussed.

One might speculate, in light of  Theorem \EE.4 above, that for minimal subequations 
the Dirichlet Problem can be solved by replacing the standard Perron family with the subfamily
consisting only of $\cp^+$-pluriharmonics. For the two extreme subequations -- the convexity and the Laplacian
subequations -- this is in fact true.  In \S 8 we prove something close for all minimal subequations.

\Theorem {\HH.3}
{\sl
Let $\cp^+$ be a minimal subequation, and $\O\ss\rn$ a domain with smooth strictly
convex boundary.  Then the standard solution to the Dirichlet problem for any $\vf\in C(\bo)$
is the upper envelope of functions in the Perron sub-family of functions  which can be written locally as 
the maximum of a finite number of pluriharmonics.
}

\vskip .2in

\medskip
\noindent{\headfont \BB.\  Preliminaries.}

In this section we review the basic properties of convex cone subequations and define
many of the associated objects (cf. [\PUP]).

 We start with a closed convex cone $\cp^+$  in $\Symn$, which we always assume 
 is a non-empty proper subset. Using the natural inner product 
 $\bra AB = \tr(AB)$ we have
 $$
 {\bf (Polar \ Cone)} \qquad\qquad \cp_+ \ \equiv\ \{A : \bra AB \geq0 \ \forall \, B\in \cp^+\}.
 \eqno{(\BB.1)}
 $$

If, in addition, $\cp^+$ satisfies the following  {positivity condition},
then $\cp^+$ is referred to as a subequation.
(We will frequently evoke the bipolar theorem, which  says that the polar of the polar is the original convex cone.)

\Def{\BB.1}  $\cp^+$ is a {\bf subequation} if satisfies the
{\bf positivity condition}
$$
\cp^+ + \cp = \cp^+, \qquad {\rm i.e.} \ \ \cp\ss \cp^+
\eqno{(P)}
$$
or equivalently   if $\cp_+\ss \cp$. The equivalence follows since $\cp$ is self-polar.

Subequations have the important {\bf topological property}
$$
\cp^+ \  = \ \overline{\Int \cp^+}
\eqno{(T)}
$$
since $\cp^+ + \e I \ \ss \ \cp^+ + \Int  \cp\  \ss \  \Int  \cp^+$ for all $\e>0$.

The following is an important class of examples.

 \Ex{\BB.2. (The Geometric Case)}  Given a closed subset $\GG\ss G(p, \rn)$
 of the Grassmannian of unoriented $p$-planes in $\rn$, let
 $$
 \cp(\GG) \ \equiv\ \{A : \bra A {P_W} = \tr(A\bigr|_W)\geq 0\ \forall\, W\in\GG\}
 $$
where $W\in \GG$ is identified with $P_W$, orthogonal projection onto $W$.
Each such $\cp^+ = \cp(\GG)$ is a  convex cone subequation with polar
$\cp_+ = $ Convex Cone Hull $(\GG)$ $\equiv {CCH}(\GG)$.

In addition to the polar cone $\cp_+$ we associate two vector spaces $E$ and $S$  with
$\cp^+$  which form an orthogonal decomposition $\Symn = E\oplus S$.
\centerline{\bf \ }
\centerline{\bf The Edge and the Span}
 $$
\text{ {\bf (The  Edge $E$)} } \ \ \ E\ \equiv \ \cp^+ \cap (-\cp^+) \ =\ \{A : A+\cp^+ = \cp^+\}
 \eqno{(\BB.2)}
 $$
 $$
\text{ {\bf (The  Dual Span $S$)}} \ \ \ \ S\ \equiv\ \span \cp_+
\qquad\qquad\qquad\qquad\qquad
 \eqno{(\BB.3)}
 $$
 Note that the edge $E$ is the unique maximal vector space contained in $\cpp$.
To verify the equality in (\BB.2), use the fact that $A+\cp^+= \cp^+ \iff -A +\cpp = \cpp$,
along with the fact that $\cpp$ is a convex cone with vertex 0.
Note that $S$ is by definition a vector space, whereas $E$ is obviously closed under addition
and scalar multiplication.

\Lemma{\BB.3}  {\sl  The Edge $E$ of a subequation enjoys the properties:

(\BB.4) \ \ (Orthogonality)  $E$ and $S$ are orthogonal compliments in $\Symn$.

(\BB.5) \ \  $E\cap (\Int\cpp) \ =\ \emptyset$.  In particular, $E\cap (\Int \cp) = \emptyset$.
}

\noindent
{\bf Proof of (\BB.4).}
It is easy to see that $E\perp S$.
Then since  $\cpm\ss S \Rightarrow S^\perp \ss \cpp$, and  
since $S^\perp$ is a vector subspace, this implies that  $S^\perp \ss E$.
Therefore, $\Symn = S^\perp +S \ss E+S$ 
thereby proving that $E+S = \Symn$ is an orthogonal decomposition.  \qed

\noindent
{\bf Proof of (\BB.5).}  If this does not hold, we can pick  $A\in E\cap \Int \cpp$.
Given $B\in \Symn$, we have $A+\e B\in \cpp$ if $\e>0$ is sufficiently small.
Therefore, $B= -{1\over \e} A +({1\over \e}A+B) \in E+\cpp =  \cpp$.  This contradicts the assumption
that $\cpp$ is a proper subset of $\Symn$.\qed
$$ 
{\rm  Let}\ \ \pi : \Symn \ \arr\ S\ \ \text{denote orthogonal projection.}
\eqno{(\BB.6)}
 $$
As a constraint on the second derivative, the important part of  $\cpp$ is
$$
\text{ {\bf (The  Reduced Constraint Set)}} \ \ \ \  \cp_0^+ \ \equiv \ \pi(\cpp).  \qquad
 \eqno{(\BB.7)}
 $$
$$
\cpp\ =\ E\oplus \cpp_0, \qquad {\rm i.e.,} \qquad A\in \cpp \quad\iff\quad \pi(A) \in \cpp_0.
 \eqno{(\BB.8a)}
 $$
 $$
\Int \cpp\ =\ E\oplus \Int \cpp_0, \quad {\rm i.e., } \quad A\in \Int \cpp \quad\iff\quad \pi(A) \in \Int \cpp_0.
 \eqno{(\BB.8b)}
 $$
Since $\pi(A)$   captures the important part of $A\equiv D^2 u$,
$\pi(D^2u)$ is called {\bf the reduced hessian of $u$ for $\cpp$.}

 Note that 
$$
\text{  $ \cp_0^+$ is the polar of $\cpm$ in its span $S$}.
 \eqno{(\BB.9)}
 $$
 
 The closed convex cone $\cp_0^+ \ss S$ is not a subequation unless  $S\equiv \Symn$, i.e., $E=\{0\}$, 
 in which case $\cp_0^+ = \cp^+$.

We say that  the subequation $\cpp$ is {\bf self polar} if
$$
\cp^+_0\ =\ \cp_+
 \eqno{(\BB.10)}
 $$

 \vskip .3in


\noindent
{\headfont \CC.  Extremely Degenerate Versus Complete}

Our subequations divide into two kinds.  The first is that of extreme degeneracy.
These subequations on  $n$-dimensional euclidean space $\rn$ are better understood
as subequations on a lower dimensional subspace (see Prop. \CC.4 below and the Support Theorem \DD.3).

\Def{\CC.1$^*$}  A convex cone subequation $\cpp$ is said to be {\bf  extremely degenerate}
if there exists a proper subspace $W\ss\rn$ such that the following equivalent conditions are satisfied:

\noindent
(1$^*$) \ The reduced constraint set  $\cpp_0 \ss \Sym(W)$,

\noindent
(2$^*$) \  The polar $\cpm \ss \Sym(W)$, or equivalently the dual span $S\ss \Sym(W)$.

\noindent
{\bf Proof that (1$^*$)$\iff$(2$^*$).}
Note that: $\cp_+ \ss\Sym(W) \iff S\equiv \span \cp_+ \ss \Sym(W) \iff \cp_0^+ \ss \Sym(W)$ by 
(\BB.10).  \qed

\noindent
{\bf Remark.}  In [\PUP] we  said for (2$^*$) that ``$\cp_+$ only involves the variables in $W$'', and for (1$^*$) that 
``$\cpp$ can be defined using the variables in $W$''.

The remaining subequations are defined by taking the negations of (1$^*$) and (2$^*$).

\Def{\CC.1}  A subequation which is not extremely degenerate will be called
{\bf dimensionally complete}, or just  {\bf complete}.  In other words

(1) $\cp_0^+ \not\ss \Sym(W)\ \ \  \text{for any proper subspace $W\ss\rn$}$, or

(2) $\cp_+ \not\ss \Sym(W)\ \ \  \text{for any proper subspace $W\ss\rn$}$.

See [\PUP] for many interesting results for complete convex cone subequations.  
The purpose  of this paper is to investigate a special class  of such subequations described in Section \DDD.
\bigskip

\centerline{\bf Extreme Degeneracy}

Let $A\bigr|_W$  denote the restriction of $A$ to the subspace $W\ss \rn$ as a quadratic form.
In terms of the $2\times 2$-blocking of $\Symn$ induced by $\rn = W\oplus W^\perp$,
$A\bigr|_W$ is the $(1,1)$-component of $A$.

The subequation $\cp^+$ can be restricted to a subequation $\cp^+_W$  on $W$
by defining
$$
\cp^+_W \ \equiv \ \left\{ A\bigr|_W : A\in \cpp\right\}.
\eqno{(\CC.1)}
$$
Note that $\cp^+_W \ss\Sym(W)$ satisfies positivity, since  if $P\in \Sym(W),$ with $P\geq0$,
then $A\bigr|_W +P = (A+Q)\bigr|_W$, where $Q\in\Symn$ restricts to $P$ on $W$
and has all other components 0, and therefore $A+Q\in \cp^+$.
Thus $\cp^+_W$ is a subequation on $W$.  

We have proved the following.
Let  $\Sym(W)^\perp$ denote the orthogonal complement of $\Sym(W)$ in $\Symn$.

\Prop{\CC.2}  {\sl
If $\cpp$ is extremely degenerate, i.e., if $\cpp_W$ and $W$ satisfy the equivalent conditions 
(1$^*$) and (2$^*$) above, then 
$$
\cpp\ =\ \cpp_W \oplus \Sym(W)^\perp.
 \eqno{(\CC.2)}
 $$
Moreover,  $\cpp$ and $\cpp_W$ have the same reduced constraint set
$\cpp_0$ since $\Sym(W)^\perp \ss E$, i.e., $S\ss \Sym(W)$. 
}

\Def{\CC.3}  If (\CC.2) is satisfied, we say that  $\cpp$ {\bf reduces to $\cpp_W$}, and that 
$\cpp$ {\bf is the trivial extension of $\cpp_W$ from $W$ to $\rn$.}

\Prop{\CC.4.\ [\DDd, Thm.\ A.4]} {\sl
Suppose that $\cpp$ reduces to a subequation $\cpp_W$ on $W$.
If $z=(x,y)\in W\oplus W^\perp=\rn$ denotes coordinates, 
 then $u(x,y)$ is $\cpp$-subharmonic if and only if for each $y$, $u(x,y)$ is
$\cpp_W$-subharmonic in $x$, but otherwise $u$ is just upper semi-continuous in $(x,y)$, i.e., 
there is no constraint on $u$ with respect to the $y$-variable.
}

\bigskip

\centerline{\bf The Edge  Criteria}

Extreme degeneracy and completeness can be described in a very simple way
in terms of edges.

\Prop{\CC.5. (The Edge Criteria)}  {\sl
The following conditions on a convex cone subequation $\cpp$ are equivalent.

(1)$\cong$(2)  \ \ $\cpp$ is complete.

(3) \ \ $E\cap \cp=\{0\}$.

(4a) \ \ $P_e \notin E$ for all $|e|=1$.

(4b) \ \ $-P_e\notin \cpp$  for all $|e|=1$.

\noindent
Stated as the edge criteria for extreme degeneracy, we have
that the following are equivalent.

(1$^*$)$\cong$(2$^*$)  \ \ $\cpp$ is extremely degenerate.

(3$^*$) \ \ $E\cap \cp \neq\{0\}$.

(4$^*$a) \ \ $P_e \in E$ for some $|e|=1$.

(4$^*$b) \ \ $-P_e\in \cpp$  for some $|e|=1$.
}

\noindent
{\bf Proof.}
We will prove the extreme degeneracy version.  First we note that 
(4$^*$a) and (4$^*$b) are equivalent.  One key to the proof is the following 
Lemma taken from [\PUP].

\Lemma{\CC.6}  {\sl
Suppose that $W$ is a hyperplane in $\rn$ with unit normal $e$.  Then}
$$
P_e\in E \qquad\iff\qquad \cp_+ \ss \Sym(W).
$$

\noindent
{\bf Proof.}

\Cor{\CC.7}  {\sl  We have that (2$^*$) $\iff$ (4$^*$a).}

\noindent
{\bf Proof.}  The only thing to note is that if condition (2$^*$), that $\cp_+\ss\Sym(W)$ for some
proper subspace $W\ss\rn$, holds,  then $\cp_+ \ss \Sym(W')$ for any hyperplane $W' \supset W$.\qed

Since (4$^*$a) implies (3$^*$) is trivial, the only thing left to prove is that (3$^*$) implies 
(4$^*$a).

\Lemma {\CC.8}  {\sl
Suppose $P\geq 0$ has null space $N\ss\rn$.  Then }
$$
P\in E \qquad\Rightarrow \qquad \Sym(N^\perp)\ \ss\ E
 \eqno{(\CC.3)}
 $$

\noindent
{\bf Proof.}  The proof is modeled on the proof of (\BB.5).
It suffices to show that  $\Sym(N^\perp) \ss\cpp$.
Given $A\in \Sym(N^\perp)$, we write $A = -tP + (A+tP)$ and note that
since $P\in E$, we have $-tP\in E$ for all $t\geq0$.
Now $A+tP\in \cp$ if $t>>0$ is sufficiently large since 
$P\bigr|_{N^\perp }$ is positive definite.  This proves that 
$A\in E+\cp\ss\cpp$.  \qed

\Cor{\CC.9}  {\sl  We have that (3$^*$) $\Rightarrow$ (4$^*$a).}

\noindent
{\bf Proof.}   If  (3$^*$) holds, choose $P\in E\cap \cp$ with $P\neq 0$.  Since $P\neq 0$,
the subspace $N^\perp\neq \{0\}$.  Pick $e\in N^\perp$ with $|e|=1$.  Then $P_e\in E$.\qed

\bigskip

\centerline{\bf The Span  Criteria}

The span criteria for completeness also provides an easy check in examples
for completeness, and has important consequences for the subequation.

\Prop{\CC.10} {\sl
The following conditions on a convex cone subequation $\cpp$ are equivalent.

(1)$\cong$(2)  \ \ $\cpp$ is complete.

(5) \ \ $S\cap (\Int \cp) \neq \emptyset$.

(6) \ \ $\Int_S \cp_+ \ss \Int \cp$.
}

\noindent
{\bf Proof.}  First we show that for a pair of vector spaces $E$ and $S$ which 
are orthogonal complements in $\Symn$, (3) and (5) are equivalent.

\Lemma {\CC.11}  {\sl
Suppose subspaces $E$ and $S$ of $\Symn$ are orthogonal  complements.  Then \smallskip

\centerline
{
 $E$ satisfies the Edge Criteria (3)
  \qquad $\iff$ \qquad  
  $S$ satisfies the Span Criteria (5).
}
}

\pf  (3) $\Rightarrow$ (5).  
If (5) is false, i.e., $S\cap  (\Int \cp) =  \emptyset$, then by the Hahn-Banach Theorem there exists
an open half-space $U$ with $S\ss \partial U$ and $\Int \cp \ss U$.  Let $N\in U$ denote the unit
normal to the hyperplane  $\partial U$.  Then $S\ss \partial U\ \Rightarrow\  N\in E=S^\perp$, while
$\Int \cp \ss U \ \Rightarrow\  \bra N P > 0 \ \forall\, P>0$,
which implies $N\in \cp$.  Thus we have $N\in E\cap \cp$, but $N\neq 0$ so that (3) is false.

(5) $\Rightarrow$ (3).  
By (5) we can pick $P\in S\cap (\Int \cp)$.  
If $A\in E\cap \cp$, then $\bra A P =0$ since $A\in E$ and $P\in S$.
However, since $A\geq 0$ and $P>0$, this implies $A=0$.\qed

\noindent
{\bf Proof that (5) $\Rightarrow$ (6).}
By (5) we can choose $P\in S\cap (\Int \cp)$.
Given $A\in \Int_{S} \cpm$, for $\e>0$ sufficiently small we have $A-\e P \in \cpm$.
Thus for all non-zero $Q \in \cp\ss\cpp$ we have $0 \leq \bra {A-\e I}{Q} = \bra AQ -\e \bra PQ$.
Since $P>0$, one has $\bra PQ > 0$, which proves that $\bra AQ > 0$ for all  non-zero $Q\geq 0$.
Thus $A>0$, which proves (6). \qed

\noindent
{\bf Proof that (6) $\Rightarrow$ (5).}  Now $\cp_+$ is a closed
convex cone in $S$.  Hence $\Int_S \cp_+ \neq \emptyset$ is equivalent to $S$ equaling the 
span of $\cp_+$, which it does by the definition of $S$.  Now pick
$P\in \Int_S \cp_+$.  Then $P\in S$ and by (6) we have $P>0$, which proves (5).\qed

This completes the proof of Proposition \CC.9.\qed

The edge and span criteria (3) and (5) for completeness motivates the following definition,
which will be used in the next section.

\Def{\CC.12} 

\noindent
 (a) A subspace $E\ss \Symn$ is called a {\bf basic edge subspace} if

(3) \ \ $E\cap \cp \ =\ \{0\}$.

\noindent
(b)  A subspace $S\ss \Symn$ is called a {\bf basic span subspace} if

(5) \ \ $S\cap (\Int \cp)\ \neq \ \emptyset$.

\noindent
(c) \ \ If in addition $E$ and $S$ are orthogonal complements, then $E, S$ well be referred to as a 
{\bf basic edge-span pair}.

\vskip.3in

\noindent
{\headfont \DD. The Supporting Subequation}

For each subequation there is a smallest subspace $W$ of $\rn$ to which the subequation reduces.

\Def{\DD.1. (Support)}  Given a convex cone subequation $\cpp$ we define the 
{\bf support of $\cpp$} to be the subspace $W\ss\rn$ which is the intersection of all subspaces
$W'\ss \rn$ which that
$$
\cpp \ =\ \cpp_{W'} \oplus \Sym(W')^\perp.
\eqno{(\DD.1)}
$$

\Lemma{\DD.2}  {\sl
The orthogonal complement of the support $W$ of $\cpp$ equals:}
$$
V \ \equiv\ \span\{e\in \rn  :  P_e\in E, |e|=1\}.
\eqno{(\DD.2)}
$$
 
 \noindent
 {\bf Proof.}   Note that (\DD.1) holds $\iff  \Sym(W')^\perp \ss E \iff S\ss \Sym(W') 
 \iff \cp_+\ss \Sym(W') \iff P_e\in E$ for all $e\perp W'$ with $|e|=1$.\qed

The support illuminates the structure of the subequation.

\Theorem{\DD.3.  (Structure Theorem)}  {\sl
Suppose $\cpp\ss \Symn$ is a convex cone subequation with support $W\ss\rn$.
Then}
$$
\cpp \ =\ \cpp_W \oplus \Sym(W)^\perp  \qquad  {\rm and}
\eqno{(\DD.3)}
$$
$$
\cpp_W \ \ss\ \Sym(W) \ \ \text{is a complete subequation.}
\eqno{(\DD.4)}
$$

 \noindent
 {\bf Proof.}  To be done later.
 
 \Def{\DD.4}  If $W$ is the support of $\cpp$, the subequation $\cpp_W$ will be called the 
 {\bf supporting subequation of $\cpp$}, and its edge $E_W$ will be called the 
{\bf supporting edge of $\cpp$}

Note that the edge of $\cpp$,
$$
E \ =\ E_W \oplus \Sym(W)^\perp,
\eqno{(\DD.5)}
$$
is larger than its supporting edge $E_W$ unless $\cpp$ is complete.

Note also that the original   subequation $\cpp$ and the supporting subequation $\cpp_W$
have the same span $S$ and the same reduced constraint set $\cpp_0$.

\vskip.3in

\noindent
{\headfont \DDD. Minimal Subequations}

These subequations are the focus of this paper.   They are all constructed as follows,
starting with a basic edge-span  pair.

\Lemma{\DDD.1} {\sl
Suppose $E, S\ss\Symn$ are orthogonal complements with $E\cap \cp = \{0\}$,  or equivalently
$S\cap (\Int \cp) \neq \emptyset$.  That is, $E,S$ is a basic edge-span pair.  Then
$$
\cpp \ \equiv\ E+\cp\ \ \text{is a subequation, and it has edge $E$ and span $S$.}
\eqno{(\DDD.1)}
$$
Moreover, if $Q^+$ is any subequation with edge $E$, then $\cpp\ss Q^+$.}

\noindent
{\bf Proof.}
Obviously $\cpp$ satisfies positivity. 
It remains to show that $\cp^+ \equiv E+ \cp$ is closed.
Let $\pi:\Symn\to S$ denote orthogonal projection as in (\BB.5).
Since $E+\cp = E\oplus \pi(\cp)$, 
$$
\text{$\cpp$ is closed if and only if $\pi(\cp)$ is closed.}
\eqno{(\DDD.2)}
$$
Now we prove that:
$$
\text{$\pi(\cp)$ is closed}.
$$
Let $K\equiv \cp\cap\{\tr = 1\}$, a compact base for $\cp$.
The image $\pi(K)$ is a compact subset of $S$.  
The basic edge condition $E\cap \cp=\{0\}$ is equivalent to 
 $0\notin \pi(K)$. This is enough to conclude that  the cone on the compact convex set $\pi(K)$ is closed.
 Thus, $\cpp \equiv E+\cp$ is a subequation. \qed

To prove that $\cpp$ has edge $E$ we must show  that
$$
\cpp \cap(-\cpp) \ =\ E \qquad{\rm or \ equivalently}\qquad \pi (\cp) \cap(- \pi(\cp)) \ =\ \{0\}.
$$
Suppose $A\in \pi (\cp) \cap(- \pi(\cp))$, i.e., $A=\pi(P_1) = -\pi(P_2)$ with $P_1, P_2\in\cp$.
Then $\pi(P_1+P_2)=0$, i.e., $P_1+P_2\in E$.   Since $E\cap \cp=\{0\}$,
$P_1+P_2=0$. But this implies 
$P_1 = P_2 =0$  and hence $A=0$. Since $\cpp$ has edge $E$, it has span $S=E^\perp$.
Finally, $\cpp\ss Q^+$, since $E\ss Q^+$ and positivity for $Q^+$ implies $\cpp \equiv E+\cp \ss Q^+$.\qed

\Def{\DDD.2}  The subequation $\cpp = E+\cp$ constructed in Lemma \DDD.1 will be 
referred to as a  {\bf minimal subequation}, or the  {\bf minimal subequation with edge $E$}.

\Cor{\DDD.3} {\sl
Suppose $\cpp$ is a minimal subequation with edge-span $E, S$. 
Then }
$$
(a)\ \ E\cap \cp \ =\ \{0\},  \qquad
(b) \ \ S\cap (\Int \cp) \ \neq\ \emptyset,  \qquad
(c)  \ \ \cpp \ \ \text{is complete}.
$$

\noindent
{\bf Proof.}
By definition of minimal we have $\cpp=E'+\cp$ where 
$E'$ satisfies (a).  By Lemma \DDD.1 the edge $E$ of $\cpp$ equals $E'$.
Lemma \CC.10 says that (a) and )b) are equivalent.
Either the edge criteria (a) $\Rightarrow$ (c), or the 
span criteria    (b) $\Rightarrow$ (c), completes the proof.\qed

There are many additional interesting properties of minimal subequations, 
besides the various completeness criteria in Section \CC.

\Theorem{\DDD.4.  (Minimality Properties)}  {\sl
Suppose $\cpp \equiv E+\cp$ is the minimal subequation  with edge $E$
and span $S$.  Then

\noindent
$(1) \ \  \cpp = E+\cp, \qquad \ \ (1a) \ \ \cpp_0 = \pi(\cp),   \qquad \ \ (1b) \ \ \cpp = E\oplus \pi(\cp)$

\noindent
$(2) \ \  \Int \cpp =   E+  \Int \cp, \ \ (2a) \ \ \Int \cpp_0 = \pi(\Int \cp), 
  \ \ (2b) \ \ \Int \cpp = E\oplus \Int \pi(\cp)$

\noindent
$(3) \ \ \cp_+ = S\cap \cp, \qquad {\rm and }\qquad (3^*) \ \    \Int_S \cp_+ =  S\cap (  \Int \cp).$
}

In fact, for complete subequations each of these eight properties characterizes minimality.

\Theorem{\DDD.5.  (Minimality Criteria)}  {\sl
Suppose $\cpp\ss \Symn$ is a complete convex cone subequation, with edge $E$
 span $S$, reduced constraint set $\cp_0$, and polar cone $\cp_+$.
 Then $\cpp$ is the minimal subequation with edge $E$
 if and only if  any one of the eight equivalent conditions in Theorem \DDD.4 hold.
}

\noindent
{\bf Proof of Theorem \DDD.4.}  Assertion (1) is by Definition \DD.2.
Next we show the following.
$$
\text{(1), (1a) and (1b) are equivalent for any subequation $\cp^+$ with edge $E$.
}
\eqno{(\DDD.3)}
$$

\noindent
(1) $\Rightarrow$ (1a):  By definition $\cpp_0 =\pi(\cpp)$. Since $\pi(E)=\{0\}$,
(1) implies that  $\pi(\cpp) = \pi(\cp)$.

\noindent
(1a) $\Rightarrow$ (1b):   This follows because $\cpp = E\oplus \pi(\cpp)$.

\noindent
(1b) $\Rightarrow$ (1):  This is obvious.

\noindent
{\bf Proof of (2).}  Obviously the open set $E+\Int \cp \ss \Int\cpp$.
If $A\in \Int \cpp$, then for small $\e>0$, $A-\e I \in \Int \cpp\ss\cpp$.  Hence there exist
$B_0\in E$ and $P\geq0$ such that  $A-\e I = B_0 + P$.  Therefore,
$A= B_0 + (P+\e I) \in E+\Int \cp$, proving that $\Int \cpp = E+\Int\cp$.

Just as in (\DDD.3), we have
$$
\text{(2), (2a) and (2b) are equivalent for any subequation $\cp^+$ with edge $E$.
}
\eqno{(\DDD.4)}
$$

\noindent
{\bf Proof of (3).} Since $0\in \cpp$ and $\cpp$ is $\cp$-monotone, we have $\cp\ss\cpp$.
Since $\cp$ is self polar, taking polars implies that $\cp_+\ss\cp$ and therefore $\cp_+\ss S\cap \cp$.

Suppose $B\in S\cap\cp$. To show $B\in \cp_+$ it suffices to show that $\bra AB\geq0$ for all $A\in\cpp$.
By minimality, if $A\in \cpp$, then $A=A_0+P$ with $A_0\in E$ and $P\in \cp$.  Now 
$\bra AB = \bra PB \geq 0$ since $\bra {A_0} B=0$.\qed

\noindent
{\bf Proof of (3$^*$).} Note that $S\cap (\Int \cp)$ is an open set in $S$, and it is contained in 
$S\cap \cp$, which is a subset of $\cp_+$ by (3).  Hence, $S\cap (\Int \cp)\ss \Int_S\cp_+$.

The only non-trivial part (and the most important part) of showing $\Int_S \cpm = S\cap (\Int \cp)$
is to show that:
$$
\Int_{S} \cpm \ \ss\ \Int \cp.
\eqno{(\DDD.5)}
$$
Since $S$ is a basic span subspace, $S\cap (\Int \cp) \neq 0$.  Choose $P\in S\cap (\Int \cp)$.
Given $A\in \Int_{S} \cpm$, for $\e>0$ sufficiently small we have $A-\e P \in \cpm$.
Thus for all non-zero $Q \in \cp\ss\cpp$ we have $0 \leq \bra {A-\e I}{Q} = \bra AQ -\e \bra PQ$.
Since $P>0$, one has $\bra PQ > 0$, which proves that $\bra AQ > 0$ for all  non-zero $Q\geq 0$.
Thus $A>0$.  \qed

\noindent
{\bf Proof of Theorem \DDD.5.} 
By Theorem \DDD.4, if $\cpp$ is minimal, then $\cpp$ satisfies each of the eight conditions.
For the converses we use the hypothesis that $\cpp$ is complete.
By the edge criteria, Proposition \CC.5(3), the edge $E$ of $\cpp$ satisfies $E\cap \cp=\{0\}$.
Therefore we can apply the construction in Lemma \DD.1 to yield a minimal subequation
$Q^+ \equiv E+\cp$ satisfying all the eight conditions.  If $\cpp$ satisfies (1)
then $\cpp = Q^+$ and so it is minimal.  Similarly, if $\cpp$ satisfies (2), then $\Int\cpp = \Int Q^+$,
so that $\cp =Q^+$ is minimal.

By (\DDD.3) we have that (1), (1a) and (1b) are equivalent.

By (\DDD.4) we have that (2), (2a) and (2b) are equivalent.

Finally, if $\cpp$ satisfies (3$^*$), then since $Q_+ = S\cap\cp$ also, we have 
$\cp_+=Q_+$ and hence $\cpp = Q^+$ is minimal.  As noted above,
 (3) $\Rightarrow$ (3$^*$).  \qed

\Remark{\DDD.6}  The property (\DDD.5) is extremely important and  useful.
See [\PUP] for more details of the following.

Given $A\geq 0$ define $\D_A u \equiv \bra {D^2 u}A$, or equivalently, from the subequation
point of view,
$$
\D_A \ \equiv\  \{B\in \Symn : \bra BA\geq 0\}.
$$
Then $u$ is $\cpp$-subharmonic if and only if $u$ is $\D_A$-subharmonic for all $A\in \Int_{\rm rel}\cp_+$.
If (\DDD.5) is true, then each such operator $\D_A$ is just a linear coordinate change of the standard
Laplacian on $\rn$ (or said differently,  it is the Laplacian on $\rn$ with a different metric).
Thus results of standard potential theory, such as $u\in L^1_{\rm loc}$, are valid for 
 $\cpp$-subharmonic functions.

One final property of minimal subequation is the following.

\Prop{\DDD.7} {\sl
Suppose $\cpp$ is a minimal subequation.  Then $\cpp$ is contained in its dual subequation}
$$
\wt{\cpp} \ \equiv\ \sim (-\Int \cpp) \ =\ -(\sim \Int\cpp).
\eqno{(\DDD.6)}
$$

\noindent
{\bf Proof.}  Since $\cpp = E+\cp$ and $\cpt^+ +\cp = \cpt^+$, it suffices to show that 
$E\equiv \cpp\cap (-\cpp) \ss\cpt^+$.
Suppose $A\notin \cpt^+$, i.e., $-A\in \Int \cpp$.  
Then by \DDD.4(2)  we have $-A = B_1+P$
with $B_1\in E$ and $P>0$.  If $A\in \cpp$ also, then $A=B_2+Q$ with $B_2\in E$ and $Q\geq0$.
Therefore, $P+Q = -B_1-B_2\in E$.  However, $P+Q>0$ contradicting Corollary \DDD.3(a).
\qed

\vskip.3in

\noindent
{\headfont \EE. Edge Functions -- Pluriharmonics}

Suppose as before that $\cpp$ is a complete convex cone subequation with edge $E$.

\Def{\EE.1}  An {\bf edge function}, or {\bf $\cpp$-pluriharmonic function}
 is  a function $u$ such that both
 $$
 \text{$u$ and $-u$ are $\cpp$-subharmonic.}
 \eqno{(\EE.1)}
 $$
Thus, by definition, $u$ is continuous.

\Def{\EE.2}  An upper semi-continuous function $u$ is  {\bf  ``sub''}  the edge functions on an
open set $X\ss\rn$ if for all domains $\O\ss\ss \rn$ and all edge functions $h$ on $\O$
which are continuous on $\ob$,
 $$
u\ \leq \ h \quad {\rm on}\ \bo 
\qquad\Rightarrow\qquad
u\ \leq \ h \quad {\rm on}\ \ob. 
 \eqno{(\EE.2)}
 $$

\Prop{\EE.3}  {\sl If $u$ is dually $\cpp$-subharmonic  on $X$, i.e., $u$ is $\cpt^+$-subharmonic for the 
dual subequation $\cpt^+$ (see (\DDD.6)), then $u$ is ``sub'' the edge functions on $X$.}

\pf  Suppose $u$ is $\cpt^+$-subharmonic
and $h$ is an edge function.  Then  $-h$ is $\cpp$-subharmonic and (\EE.2) follows from comparison
(see Thm. 6.2 in [\AE]).\qed

Now if a  subequation becomes smaller, it dual subequation becomes larger.
Consequently, the only subequation $\cpp$, with a given edge $E$, for which Proposition \EE.3
might have a converse is the minimal subequation with edge $E$ (see Definition \DDD.2).

\Theorem{\EE.4}
{\sl
Suppose that $E\ss \Symn$ is a basic vector subspace, 
so that $\cpp\equiv E+\cp$ is the minimal subequation with edge $E$. Then the following conditions on a function
$u$ are equivalent.
$$
\begin{aligned}   
&\text{(1)\ \ $u$ is dually $\cpp$-subharmonic}. \\ 
&\text{(2)\ \ $u$ is ``sub'' the edge functions}. \\
 &\text{(3)\ \ $u$ is locally ``sub'' the edge functions}. \\
 &\text{(4)\ \ $u$ is locally ``sub'' the degree-2 polynomial edge functions}.
\end{aligned}
$$
}

\pf Because of Proposition \EE.3 we need only prove that if $u$ is locally  ``sub'' the degree-2 polynomial edge 
functions, then $u$ is dually $\cpp$-subharmonic.
For this suppose  that $u$ is not $\cpt^+$-subharmonic on $X$.
Then (see Lemma 2.4 in  [\DDR]) there exists $z_0\in X$, a quadratic polynomial test  function  $\vf$, and $\a>0$
such that
$$
u(z) \ \leq \ \vf(z) -\a|z-z_0|^2\quad {\rm near}\ \ z_0\ \ {\rm with\ equality\ at\ \ } z_0,
\eqno{(\EE.3)}
$$
but
$$
D^2_{z_0} \vf \notin \cpt^+, \ \ \  {\rm i.e., }\ \ -D^2_{z_0}\vf \in   \Int \cpp.
\eqno{(\EE.4)}
$$
By Theorem \DD.4(2) we have $\Int \cpp = \Int \cp+E$. Thus
$$
 -D^2_{z_0}\vf \ =\ P+B  \ \ \ {\rm with}\ \ \ P>0\ \ \ {\rm and}\ \ \ B\in E.
\eqno{(\EE.5)}
$$
Consider the  degree-2  edge polynomial   
$$
\begin{aligned}
h(z) \  &\equiv\  \vf(z_0) + \bra {D_{z_0} \vf}{z-z_0}  -  \half\bra {B(z-z_0)}{z-z_0}   \cr
&= \ \vf(z) -\half \bra  {D_{z_0}^2 \vf}{z-z_0}  -  \half\bra {B(z-z_0)}{z-z_0}   \cr
&=\  \vf(z) + \half \bra {P(z-z_0)}{z-z_0}.
\end{aligned}
$$
Since $P>0$ by (\EE.3) this implies that 
$$
u(z) \ \leq\ h(z) - \a|z-z_0|^2
\eqno{(\EE.6)}
$$
near $z_0$ with equality at $z_0$.  This implies that  $u$ is not sub the function $h$ 
on any small ball about $z_0$.  Hence, $u$ is not locally ``sub'' the  degree-2 edge  polynomial $h$.
\qed

\vskip.3in

\noindent
{\headfont \FF. Further Discussion of Examples}

Before turning to the examples we define the {\bf (compact) invariance group of $\cpp$} to be
$$
\{ g\in {\rm O}_n : g^*\cpp \ =\ \cpp\}.
\eqno{(\FF.1)}
$$
  It is easy to see that for the minimal subequation $\cpp$ for a basic $E$, 
$$
g^*\cpp= \cpp  
\iff
g^*S=S
 \ \iff\ 
g^* E= E
\eqno{(\FF.2)}
$$
by using the conditions in Theorems \DDD.4 and \DDD.5,
and this yields two equivalent definitions of this group.

\Def{\FF.1. (Self Duality)}  If the two convex cones $\cpp_0$ (the reduced constraint set) and 
$\cpm$  are polars of each other  in the vector space $S$, then we say the subequation $\cpp$ is {\bf polar self dual}
(not to be confused with a subequation which equals its dual subequation in the sense of [\DDd]).

\Remark  {\FF.2}  Note that this can only happen for a minimal subequation $\cpp$.
This is because if $\cpp_0 = \cpm$ (self duality), then $\cpp_0= \cpm \ss \cp$,
and hence $\cpp_0 = \pi(\cpp_0) \ss \pi(\cp)$.
Note that $\cp\ss \cpp$ so that $\pi(\cp) \ss\cpp_0$ is always true.
This proves $\cpp_0 = \pi (\cp)$, so by Theorems \DDD.5 and \DDD.4(1a), $\cp^+$ is minimal.

Given a closed subset $\GG\ss G(k,\bbr^n)$ consider  the {\sl subequation geometrically defined by $\GG$}:
$$
\cp(\GG) \ \equiv\ \left \{ A\in \Sym(\bbr^n) : \bra A {P_W} = \tr (A\bigr|_W)\geq 0 \ \ \forall\, W\in\GG \right\}.
$$
We shall use the following notations introduced in Example \BB.2:
$$
\cp^+ \equiv \cp(\GG), \ \ \cp_+ = CCH(\GG),\ \ S= \span(\GG),\ \ E= S^\perp, \ \ {\rm and}\ \ \cpp_0.
$$
 Note that the compact invariance group of the subequation $\cpp = \cp(\GG)$ can also be
defined by
$$
\{ g\in {\rm O}(n) : g(\GG)=\GG\}.
\eqno{(\FF.2)'}
$$

\vskip .2in

\centerline
{
\bf  The O(n)-Invariance Group
}

For our first two examples of  minimal  subequations we focus on the 
${\rm O}_n$-{\bf orthogonal decomposition}
$$
\Symn \ =\ \bbr \cdot {\rm Id} \oplus \Sym_0(\rn)
\qquad\qquad\qquad 
\eqno{(\FF.3)}
$$
into irreducible components under O$_n$.

\Ex{\FF.1.  (Real Monge-Amp\`ere)} The subequation is $\cpp=\cp$.  Here the edge $E=\{0\}$ is as small
as possible, and $S=\Symn$, $\cpm=\cp$, so the subequation is self-dual, and we have $\GG=G(1, \rn)$.
Obviously $E, S$ is a basic edge-span pair (Definition \CC.12c). 
The conditions in Theorem \DDD.4  are obvious as well as  the fact that 
$\cp=\cpp=\cpm=\cpp_0$ is dimensionally complete. The invariance group is O$_n$, and the extreme rays
are
$$
\Ext(\cp) \ =\ \{{\rm Ray}(P_e) : |e|=1\}.
$$
 Each $A\in S$ can be put in canonical form $A=\sum_j \la_j P_{e_j}$ under the action of O$_n$, 
 and $\det (A) = \prod_j \la_j$, 
  provides a nonlinear operator for $\cpp =\{ \la_{\rm min}\geq 0\}$ 
 (the standard real Monge-Amp\`ere operator).

\Ex{\FF.2.  (The Laplacian)}  Here $\cpp= \D =\{A : \tr(A)\geq0\}$ is a closed half space, and
$\GG=\{{\rm Id}\} = G(n, \rn)$,  $E = \Sym_0(\rn)$, the traceless part of $\Symn$, $S=\bbr\cdot \Id$,
and $\cpm= \bbr_+\cdot \id$ is a ray.  The invariance group is O$_n$.  The reduced constraint set 
is $\cpp_0 = \cpm$ so $\D$ is self dual.  Now it is obvious that $\D$ is a  minimal subequation.

\vskip .2in

\centerline
{
\bf  The U(n)-Invariance Group
}

We now consider $\cn$ and the following U$(n)$-{\bf orthogonal decomposition}
 of real symmetric matrices into U$_n$-irreducible subspaces:
$$
\Sym_\bbr(\cn) \ =\ \bbr\cdot\Id \oplus \Herm_0^{\bbc-{\rm sym}}(\cn) 
\oplus \Herm^{\bbc-{\rm skew}}(\cn) 
\quad
\eqno{(\FF.4)}
$$
multiples of the identity, traceless complex hermitian symmetric, and complex hermitian skew components. 
Given $A\in \Sym_\bbr(\cn)$, this decomposition can be written as
$$
A\ =\ {\tr(A) \over 2n} \Id + A_0^{\bbc-{\rm sym}}  + A^{\bbc-{\rm skew}}
\eqno{(\FF.5)}
$$
where with respect to multiplication $I$ by $i$:
$$
 A^{\bbc-{\rm sym}} \ =\ \half (A-IAI) \and 
 A^{\bbc-{\rm skew}}\ =\ \half (A+IAI).
$$

\Ex{\FF.3. (Complex Plurisubharmonics)}  The subequation is $\cpp= \cp(\GG)$
where $\GG = \bbp(\cn) \ss G_\bbr(2, \cn)$ is the Grassmannian of complex lines in $\cn$.
The edge is $E=  \Herm^{\bbc-{\rm skew}}(\cn) $ and the span is $S = \Herm^{\bbc-{\rm sym}}(\cn)$.
Also $\cpp_0=\cpm$ is the convex cone on non-negative complex hermitian symmetric bilinear
forms on $\cn$,  so this third example is  self dual.
  Note that the projection of $2 P_e$ onto $S$ is $P_{\bbc e}$,  (orthogonal  projection onto the
complex line through $e$) since $P_e - IP_e I =P_{\bbc e}$.
The convex cone $\cpp_0 =\cpm$ has extreme rays generated by
$\{P_{\bbc e} : |e|=1\} = \bbp(\cn)=\GG$.  The invariance group
is U$_n$.  Each $A\in S$ can be put into canonical form $A=\sum_{j=1}^n \la_j P_{\bbc e_j}$
under the action of this group, and $\cpp_0 = \{\la_{\rm min}\geq 0\}$.
The complex Monge-Amp\`ere operator $\det(A) = \la_1(A)\cdots \la_n(A)$ provides the nonlinear operator for
$\cpp=\cp(\bbp(\cn))$, in tight analogue with the real case $\cp$.

Now we finally get to a new example, which is the subject of [\LA].

\Ex{\FF.4. (Lagrangian Plurisubharmonics)}  The subequation is $\cpp=\cp(\LAG)$,
where $\LAG\ss G_\bbr (n,\cn)$ is the set of Lagrangian $n$-planes in $\cn=\bbr^{2n}$.
The edge $E$ and span $S$ are given by
$$
E=\Herm_0^{\bbc-\sym}(\cn)
\and
S=\bbr\cdot \Id \oplus \Herm^{\bbc-\sk}(\cn).
$$
In [\LA] we prove that   $E,S$ is a basic edge-span pair, so that $\cpp=E+\cp$ and
$\cpm = S\cap \cp$.  The extreme rays in $\cpm$ are generated by the projections $P_W$
with $W\in \LAG$ a Lagrangian $n$-plane.  The extreme rays in $\cp_0^+$ are generated
by the images  $\pi(P_e)$ of $P_e$ where $e$ is a unit vector.  Note that 
$$
\pi(P_e) = \smfrac{1}{ 2n} \Id +\half (P_e+IP_e I) = \smfrac{1}{ 2n} \Id +\half (P_e - P_{Ie}),
$$
and that $\half (P_e - P_{Ie})$ is the $\bbc$-skew component of $P_e$.
This example is {\sl not} self dual.  
However, since each $A\in S$ can be put in canonical form
$$
A\ =\ {\tr(A)\over 2n} + {1\over 2} \sum_{j=1}^n \la_j \left(  P_{e_j}- P_{Ie_j}\right)
$$
there is again a nonlinear operator for $\cpp=\cp(\LAG)$ (see [\LA]).
The invariance group is U$_n$.

\vskip .2in

\centerline
{
\bf  The Sp(n)$\cdot$Sp(1)-Invariance Group
}

Let $M_n(\bbh)$ denote the space of $n\times n$ matrices with entries in $\bbh$, and let
$A^* = \overline A^t$ if $A\in  M_n(\bbh)$.  Consider the two subspaces
$$
M_n^\sym(\bbh) \ =\ \{A \in M_n(\bbh) : A^*=A\}, \qquad {\rm and}
$$
$$
M_n^\sk(\bbh) \ =\ \{A \in M_n(\bbh) : A^*= - A\}.  \qquad \ \ \ 
$$
We let the scalars $\bbh$ act on the right.  Then by letting $M_n(\bbh)$
act on $x=(x_1, ... ,x_n)^t\in \hn$ on the left, one can identify $M_n(\bbh)$ with  ${\rm End}_\bbh(\hn)$,
the vector space of $\bbh$-linear maps of $\hn$.
Let 
$$
\Herm^{\bbh-\sym}(\hn) \ = \ \{A\in {\rm End}_\bbh(\hn) : A=A^*\}, \qquad  {\rm and}
$$
$$
\Herm^{\bbh-\sk}(\hn) \ = \ \{A\in {\rm End}_\bbh(\hn) : A= - A^*\}. \qquad  {\rm \ \ \ }
$$
so that $M_n^\sym(\bbh) = \Herm^{\bbh-\sym}(\hn)$ are identified (same for the skew parts).

Let $\e(x,y) = \sum_{\ell=1}^n \overline x_\ell y_\ell$ denote the {\sl standard quaternionic hermitian bilinear form on }
$\hn$.  The {\bf quaternionic unitary group} is 
$$
{\rm Sp}_n \ =\ \{A\in M_n(\bbh) : \e(Ax, Ay) = \e(x,y)\}.
$$
For each scalar $u\in \bbh$ let $R_u x\equiv xu$ denote right multiplication, and set
$I\equiv R_i$, $J\equiv R_j$, $K\equiv R_k$.  Then the group of unit scalars 
${\rm Sp}_1 \equiv S^3 = \{R_u : u\in \bbh, |u|=1\}$ acts on $\hn$ on the right and the 
{\bf enhanced quaternionic unitary group} is the group
$$
{\rm Sp}_n \cdot {\rm Sp}_1 \ =\ {\rm Sp}_n  \times {\rm Sp}_1/\bbz_2.
$$

Since the standard euclidean inner product on $\bbr^{4n}=\hn$ is $\bra xy = {\rm Re} \,\e(x,y)$, 
$$
M_n^\sym(\bbh) = \Herm^{\bbh-\sym}(\hn) \ \ \text{is a real subspace of $\Sym(\bbr^{4n})$}
$$
and
$$
M_n^\sk(\bbh) = \Herm^{\bbh-\sk}(\hn) \ \ \text{is a real subspace of ${\rm  Skew}^2(\bbr^{4n})$}
$$
where End$_\bbr(\bbr^{4n}) = \Sym(\bbr^{4n})\oplus {\rm  Skew}^2(\bbr^{4n})$ is the usual decomposition.
Note also that for each unit imaginary quaternion $u\in {\rm Im}\bbh$, we have 
$R_u \in  {\rm  Skew}^2(\bbr^{4n})$, and hence $R_uA = A R_u \in \Sym(\bbr^{4n})$
for all $A \in M_n^\sk(\bbh) = \Herm^{\bbh-\sk}(\hn)$.  This embeds
$$
{\rm Im} \bbh \otimes \Herm^{\bbh-\sk}(\hn)  \ =\ {\rm Im} \bbh \otimes M_n^\sk(\bbh) \ \ss\ \Sym(\bbr^{4n}).
\eqno{(\FF.6)}
$$

The ${\rm Sp}_n \cdot {\rm Sp}_1$-{\bf orthogonal decomposition}
$$ 
\Sym(\bbr^{4n})  = \bbr\cdot \Id \oplus \Herm_0^{\bbh-\sym}(\hn)\oplus  
\left(  {\rm Im} \bbh\otimes \Herm^{\bbh-\sk}(\hn) \right)
\eqno{(\FF.7)}
$$
into irreducible components plays a role in the next two examples,
and a key role in classifying all the Sp$_n \cdot {\rm Sp}_1$-invariant minimal subequations.
Projection onto $\Herm^{\bbh-\sym}(\hn) =  \bbr\cdot \Id \oplus \Herm_0^{\bbh-\sym}(\hn)$
and ${\rm Im} \bbh\otimes \Herm^{\bbh-\sk}(\hn)$ are given by $A= A^{\bbh-\sym} + A^{\bbh-\sk}$ where
$$
A^{\bbh-\sym}  \ =\ \smfrac 14 (A-IAI-JAJ-KAK),\qquad{\rm and}
\eqno{(\FF.8a)}
$$
$$
A^{\bbh-\sk}  \ =\ \smfrac 14 (3A+IAI+JAJ+KAK).\qquad\ \ \ 
\eqno{(\FF.8b)}
$$

\Ex{\FF.5. (Quaternionic Plurisubharmonics)}   The subequation $\cpp\equiv \cp(\bbp(\hn))$
is geometric with $\GG= \bbp(\hn)\ss G_\bbr(4, \hn)$.  The edge and span are given by
 $$
 E= {\rm Im} \bbh \otimes \Herm^{\bbh-\sk}(\hn),
\and S= \Herm^{\bbh-\sym}(\hn).
\eqno{(\FF.9)}
$$
  The set $\cpp_0=\cpm$ is the convex cone of
non-negative quaternionic hermitian symmetric bilinear forms on $\hn$ (see [\Bry] or [\Har] for more 
details).  Under the identification of $\Herm^{\bbh-\sym}(\hn)$ with the set of quaternionic $n\times n$ matrices 
$M_n(\bbh)$ satisfying $A^* \equiv {\overline A}^t = A$, we have
$$
\cpp_0 \ =\ \{A\in M_n(\bbh) : \ A^* = A\ \ {\rm and}\ \ \overline x^t A x  \geq0 \ \forall\, x\in \hn\}.
$$
This is a minimal subequation and has compact invariance group ${\rm Sp}_n \cdot {\rm Sp}_1$.
Note that by (\FF.7a) the projection of $P_e$ ($|e|=1$) onto $\Herm^{\bbh-\sym}(\hn)$ is just
$P_{\bbh e}$, orthogonal projection onto the quaternionic line $\bbh e$.  Hence, this example is self dual,
i.e., $\cpp_0 = \cpm$.  Each $A\in M_n(\bbh)$ with $A^*=A$ can be put in canonical form
$$
A^{\bbh-\sym} \ =\ \sum_{j=1}^n \la_j P_{\bbh e_j}
$$
under the action of ${\rm Sp}_n \cdot {\rm Sp}_1$.
The quaternionic Monge-Amp\`ere operator 
$$
\det_\bbh (A) \ \equiv \ \prod_{j=1}^n \la_j \left (A^{\bbh-\sym} \right)
$$
provides the nonlinear operator  for  $\cpp = \cp(\bbp(\hn))$.

\Ex{\FF.6a}  Reversing the roles of ${\rm Im} \bbh \otimes  \Herm^{\bbh-\sk}(\hn)$
  and  
  \\ \noindent $\Herm^{\bbh-\sym}_0(\hn)$ in (\FF.9) above  results in a second
${\rm Sp}_n \cdot {\rm Sp}_1$-invariant minimal subequation $\cp^+\equiv E+\cp$
with $\cpm = S\cap \cp$, where
$$
E \ \equiv\ \Herm^{\bbh-\sym}_0(\hn) \ \ {\rm and}\ \  S \ \equiv \  \bbr\cdot \Id \oplus 
({\rm Im} \bbh \otimes  \Herm^{\bbh-\sk}(\hn)).
\eqno{(\FF.10)}
$$
Note that  for each $|e|=1$,
$$
\pi(P_e) \ =\ \smfrac 1 {4n} \Id  + \smfrac 14 (3P_e -P_{Ie} -P_{Je} -P_{Ke}).   
\eqno{(\FF.11)}
$$
We leave as a question: Does $\pi(P_e)$ generate an exposed ray in $\cp^+_0 = \pi(\cp)$?

This edge $E\equiv \Herm_0^{\bbh-\sym}(\hn)$  is reminiscent of the edge in Example \FF.4 in the complex case.  We now pursue this analogy.  We say
that a real $n$-plane $W$ in $\hn$ is  {\bf $\bbh$-Lagrangian} if 
$$
W \oplus IW\oplus JW \oplus KW \ =\ \hn \quad\text{(orthogonal direct sum)},
\eqno{(\FF.12)}
$$
and let $\bbh\Lag$ denote the set of all such $n$-planes.

\Ex{\FF.6b. (Quaternionic Lagrangian Plurisubharmonics)} 
These are defined as the subharmonics for the geometrically defined subequation $\cp(\bbh\Lag)$.
Note that $\bbh\Lag$ and hence $\cp(\bbh\Lag)$ has compact invariance group 
${\rm Sp}_n \cdot {\rm Sp}_1$.  Furthremore, given $A\in \Sym_\bbr(\hn)$ one can show that 
$$
\tr A\bigr|_W\ =\ 0 \quad \forall\, W\in \bbh\Lag 
\qquad\iff\qquad A\in \Herm^{\bbh-\sym}_0(\hn)).
\eqno{(\FF.13)}
$$
Consequently, $\cp(\bbh\Lag)$ has edge-span given by (\FF.10).  At the moment we
do not know whether or not $\cp(\bbh\Lag)$ is the minimal subequation with this 
edge-span.  Of course one has
$$
E+\cp \ \ss\ \cp(\bbh\Lag) \and \cp_+(\bbh\Lag) \ \ss \  S\cap \cp,
\eqno{(\FF.14)}
$$
where $\cp_+(\bbh\Lag)$ is the convex cone hull of $\{P_W : W\in \bbh\Lag\}$.

\vskip.3in

\centerline{\bf The Sp$_n\cdot$S$^1$--Invariance group}

If  U$_n$ is replaced by the smaller subgroup SU$_n$, the decomposition (\FF.4) of 
$\Sym_\bbr(\cn)$ remains the same, and so SU$_n$ is not a compact invariance group 
for a minimal subequation.  However, the decomposition (\FF.7) does not remain
the same if we replace ${\rm Sp}_n \cdot {\rm Sp}_1$
by Sp$_n$.  The new decomposition can be written as
$$
{\rm Sp}_n: \quad
\Sym_\bbr(\hn) \ = \   
\bbr\cdot \Id     \oplus
\Herm^{\bbh-\sym}_0(\hn)   
\bigoplus_{j=1}^3 I_j \Herm^{\bbh-\sk}(\hn)   
\eqno{(\FF.15)}
$$
where $I_j$ vary over $I,J,K$, or in fact over any orthonormal basis of ${\rm Im} \bbh$.
Note that the representations $ I_j \Herm^{\bbh-\sk}(\hn)$ are all equivalent.

The next example is a minimal subequation which is new.

\Ex{\FF.7. ($I$-Complex and $J, K$-Lagrangian Plurisubharmonics)}   This is a geometrically
defined subequation given by the set
$$
\GG \ =\ \GG (I; J,K) \ \ss\ G_\bbr(2n, \hn)
$$
of real $2n$-planes which a simultaneously $I$-complex and both $J$ and $K$ Lagrangian.
(Note that any two of these conditions implies the third.) The associated subequation is
$\cp(\GG (I; J,K))$.  

Now $\cp({\rm JLAG})$ has edge 
$$
\Herm^{J\bbc-\sym}_0(\bbc^{2n}) =  \Herm^{\bbh-\sym}_0(\hn)  \oplus J \Herm^{\bbh-\sk}(\hn)
$$
and $\cp({\rm KLAG})$ has edge 
$$
\Herm^{K\bbc-\sym}_0(\bbc^{2n}) =  \Herm^{\bbh-\sym}_0(\hn)  \oplus K \Herm^{\bbh-\sk}(\hn).
$$
Hence the sum
$$
 \Herm^{\bbh-\sym}_0(\hn)  \oplus J \Herm^{\bbh-\sk}(\hn)\oplus  K \Herm^{\bbh-\sk}(\hn)
 \ \ss\ {\rm Edge}(\cp(\GG)).
$$

Each $W\in \GG$ has a real basis of the form
$$
e_1, Ie_1, ... , e_n, Ie_n \qquad
\text{where $e_1,...,e_n$ is an $\bbh$-basis for $\hn$.}
$$
Thus $P_W = P_V +P_{IV}$ where $V\equiv \span_\bbr \{e_1, ... , e_n\}$.
Note  that $W\in \GG \Rightarrow W^\perp = JW = KW\in \GG$.  Hence, $\Id = P_W+P_{W^\perp}
\in S \equiv \span(\cp(\GG)) \equiv \span(\GG)$.
Now we have 
$$P_W - P_{W^\perp} = P_W - P_{IW} = P_W + IP_{W}I  \ \in \  I  \Herm^{\bbh-\sk}(\hn).$$
One can show (direct proof and invariance proof) that
$$
S\ = \ \bbr\cdot\Id \oplus  I \Herm^{\bbh-\sk}(\hn)
\eqno{(\FF.16)}
$$
and hence
$$
E \ =\  \Herm^{\bbh-\sym}_0(\hn)  \oplus J \Herm^{\bbh-\sk}(\hn) \oplus K \Herm^{\bbh-\sk}(\hn).
\eqno{(\FF.17)}
$$

\Lemma{\FF.8}  {\sl
Each $A\in I \Herm^{\bbh-\sk}(\hn)$ commutes with $I$ and anti-commutes with $J$ and $K$.
If $e$ is an eigenvector of $A$  with eigenvalue $\la$, then $Ie, Je, Ke$ are eigenvectors with
eigenvalues $\la, -\la, -\la$.  Hence, $A$ can be put in the canonical form (where $e_1, ... , e_n$ is an
$\bbh$-basis for $\hn$):}
$$
A \ \equiv\ \sum_{j=1}^n \la_j \left( P_{e_j} +  P_{Ie_j} -  P_{Je_j} -  P_{Ke_j}\right).
$$ 
 
\Cor{\FF.9}  {\sl
The element $B \equiv {t\over 4n}\Id +A \in S$ is $\geq0$  if and only if  each $|\la_j |\leq {t\over 2n}$.
Hence, taking $t = \tr(B) = 2n$, the non-negativity condition becomes
$$
|\la_j|\ \leq\ {1\over 2}, \quad j=1, ... , n.
$$
}

This describes a cube in $\rn$.  The $2^n$ extreme points are $\e = (\pm \half, ... , \pm \half)$, which yields
$$
B(\e) \ \equiv\ \half\Id + \sum_{j=1}^n \pm\half  \left( P_{e_j} +  P_{Ie_j} -  P_{Je_j} -  P_{Ke_j}\right) 
\ =\  P_{W(\e)}
$$
where
$$
W(\e) \ =\ \span\biggl\{ (e_1, Ie_1 \ \text{if } \ \e_1 = \half) \ \ \text{or}\ \ (Je_1, Ke_1 \ \text{if } \ \e_1 = -\half), \ \ 
... \ \ \text{etc.}
\biggr\}.
$$
This proves

\Prop{\FF.10}
$$
 S\cap \cp \ =\ { CCH}\{P_W : W\in \GG\} \ \equiv\ \cpm(\GG).
$$

\Cor{\FF.11} {\sl
The subequation 
$\cp(\GG)$ is {\bf the} minimal subequation with span $S$
 and edge $E$ given by (\FF.9) and (\FF.10).
}

\Ex{\FF.12}  
Consider the edge $E_I \equiv I\Herm^{\bbh-\sk}(\hn)$ 
and the minimal subequation $\cp^+\equiv E_I + \cp$.
The compact invariance group is Sp$_n\cdot S^1$,  as in Example \FF.7.

\Lemma{\FF.13} {\sl One has
$$
\begin{aligned}
\cpp \ \equiv\ E_I+\cp \ &\ss\ \cp  (I\Lag)\cap \cp( \bbp_J(\bbc^{2n})) \cap\cp( \bbp_K(\bbc^{2n}) )  \\
&=\  \cp  \left( (I\Lag)\cup \bbp_J(\bbc^{2n}) \cup \bbp_K(\bbc^{2n})  \right) 
\end{aligned}
$$
which has edge $E_I$.}

\pf
Suppose   for all $W\in I\Lag \cup \bbp_J(\bbc^{2n}) \cup \bbp_K(\bbc^{2n})$ that $\bra A {P_W} \geq0$.
Taking $W\in  I\Lag$ proves that $A\in\cp( I\Lag)$;
 taking $W\in  \bbp_J(\bbc^{2n})$ proves  that $A\in\cp(\bbp_J(\bbc^{2n}))$;  and 
 taking $W\in  \bbp_K(\bbc^{2n})$ proves  that $A\in\cp(\bbp_K(\bbc^{2n}))$.
Conversely, if $A$ belongs to the intersection of the three geometric subequations in
the Lemma, then $\tr A\bigr|_W\geq0$ for all $W\in I\Lag \cup \bbp_J(\bbc^{2n}) \cup \bbp_K(\bbc^{2n})$.
This proves the last equality in the Lemma.

Since $E_I \ss E_{0,I}$,  by Example \FF.4,
$$
\cpp \ \equiv\ E_I+ \cp \ss E_{0,I}+\cp \ =\ \cp(I\Lag).
$$
Since $E_I \ss E_{I,K}$,  by Example \FF.3,
$$
\cpp \ \equiv\ E_I+ \cp \ss E_{I,K}+\cp \ =\ \cp(\bbp_J(\bbc^{2n})).
$$
Since $E_I \ss E_{I,J}$,  by Example \FF.3,
$$
\cpp \ \equiv\ E_I+ \cp \ss E_{I,J}+\cp \ =\ \cp(\bbp_K(\bbc^{2n})).
$$
Finally since $E_I = E_{0,I}\cap E_{I,K}\cap E_{I,J}$, this proves that
$\cp  \left( (I\Lag)\cup \bbp_J(\bbc^{2n}) \cup \bbp_K(\bbc^{2n})  \right)$
has edge $E_I$.\qed

It remains an open question whether or not 
$\cp  \left( (I\Lag)\cup \bbp_J(\bbc^{2n}) \cup \bbp_K(\bbc^{2n})  \right)$
is the minimal subequation $\cpp\equiv E_I+\cp$ with edge $E_I$.

\vskip.3in

\noindent
{\headfont \GGG. Classifying the Invariant Minimal Subequations}

Given a compact subgroup $G\ss {\rm O}_N$,  one could ask which (if any) 
subequations have $G$ as their exact invariance group.
Now the compact invariance group for a minimal subequation $\cpp=E+\cp$
is the same as for its edge $E$ (see (\FF.2)).  Therefore we need only classify the possible
invariant edges $E$.  This is easily done as follows.  
First decompose $\Sym_\bbr(\bbr^N)$ into irreducible pieces
$\Sym_\bbr(\bbr^N) = \bbr\cdot \Id \oplus E_0\oplus E_1\oplus \cdots\oplus E_k$,
and note that $E_0\oplus  \cdots\oplus E_k = \Sym_0(\bbr^N)$, the traceless part.
Hence any space $E= E_{i_1}\oplus \cdots \oplus E_{i_\ell}$, $0\leq i_1<\cdot<\i_\ell\leq  k$
can be chosen as a basic (invariant) edge.  Note that $E=\{0\}$ is also a basic invariant edge,
and $E+\cp=\cp$, which has compact invariance group O$_N$.

\noindent
{\bf The O$_n$-Case.}  Here we have 
$$
\Symn \ =\ \bbr\cdot \Id \oplus E_0 \qquad{\rm with}\qquad E_0 \ \equiv\    \Sym_0(\bbr^n).
$$
There are two examples: $E=\{0\}$ and $E = E_0$  given by Examples \FF.1 and \FF.2.

\noindent
{\bf The U$_n$-Case.} Here it is more complicated:
$$
\Sym_\bbr(\cn)\ =\ \bbr\cdot \Id \oplus E_0 \oplus E_1  \qquad{\rm with}\qquad
$$
$$
E_0 \ \equiv \ \Herm_0^{\bbc-\sym}(\cn) \and  E_1 \ \equiv\    \Herm^{\bbc-\sk}(\cn),
$$
which are Examples \FF.3 and \FF.4.

\noindent
{\bf The Sp$_n\cdot$Sp$_1$-Case.} Here  we have
$$
\Sym_\bbr(\hn)\ =\ \bbr\cdot \Id \oplus E_0 \oplus E_1  \qquad{\rm with}\qquad
$$
$$
E_0 \ \equiv \ \Herm_0^{\bbh-\sym}(\hn) \and  E_1 \ \equiv\  {\rm Im} \bbh \otimes  \Herm^{\bbh-\sk}(\hn).
$$
Hence again there are two new examples $E=E_0$ and $E=E_1$ which are Examples \FF.5 and \FF.6a.

\noindent
{\bf The Sp$_n$ and  Sp$_n\cdot$S$^1$-Cases.} Under Sp$_n$ we have
$$
\Sym_\bbr(\hn)\ =\ \bbr\cdot \Id \oplus E_0 \oplus E_I  \oplus E_J  \oplus E_K  \qquad{\rm with}\qquad
$$
$$
E_0 \ \equiv \ \Herm_0^{\bbh-\sym}(\hn), \ \ \ 
 E_I \ \equiv\  I    \Herm^{\bbh-\sk}(\hn), \ \ \ $$ $$
 E_J \ \equiv\  J    \Herm^{\bbh-\sk}(\hn), \ \ \ 
 E_K \ \equiv\  K    \Herm^{\bbh-\sk}(\hn), \ \ \ 
$$
(see (\FF.15)).  Of the possible edges we can exclude most of them as coming from the previous cases.
For example, $E_{0,I} \equiv E_0\oplus E_I =  \Herm_0^{\bbc-\sym}(\cn)$ for the $I$-complex case  (as well as
$E_{0,J}$, $E_{0,K}$) come from Example \FF.4.  The case $E_{J,K} \equiv E_J\oplus E_K = \Herm^{\bbc-\sk}(\cn)$
 (for the complex structure $I$)  can be excluded, since this is Example \FF.3. Similarly we exclude $E_{I,K}$ and $E_{I,J}$.
The case $E=E_0$ is just Example \FF.6a, while the case $E\equiv E_{I,J,K}= E_I\oplus E_J\oplus E_K
=  {\rm Im} \bbh \otimes  \Herm^{\bbh-\sk}(\hn)$ is Example \FF.5.  

This leaves, up to permuting $I,J,K$, two examples:  $E=E_I$, which is Example \FF.12, and
$E = E_{0,J,K} = E_0\oplus E_J \oplus E_K$ as in (\FF.10),
which is Example \FF.7.  These last two examples have compact invariance group 
Sp$_n\cdot$S$^1$.  Note that this proves that there are no minimal subequations 
with compact invariance group Sp$_n$.

\vskip.3in

\noindent
{\headfont   \HH.  An Envelope Problem for Minimal Subequations.}

Suppose that $\bbf \equiv \cp^+$ is a minimal subequation.  In this section we investigate
the role played by the edge functions in solving the Dirichlet problem. 
 The key fact about $\bbf$ that will be used below is the following from Theorem \DDD.4(2):
 $$
 \Int \bbf \ =\ E + \Int \cp.
\eqno{(\HH.1)}
$$

We recall that existence and uniqueness for the (DP) on a   bounded domain $\O\ss\rn$ 
and arbitrary $\vf\in C(\bo)$ was established in [\DDd] if $\bo$ is smooth and strictly
$\bbf$- and $\wt \bbf$-convex (for any subequation $\bbf\ss\Symn$).  Moreover, the solution $H$
equals the Perron function
$$
H(x) \ \equiv\ \sup_{u\in \cf_{\bbf}(\vf)} u(x) \qquad{\rm for}\ \ x\in \ob
\eqno{(\HH.2)}
$$
for the Perron family of $\bbf$-subharmonics
$$
\cf_{\bbf}(\vf) \ \equiv\ \left \{u\in \bbf(\ob) : u\bigr|_{\bo}\leq\vf\right\}.
\eqno{(\HH.3)}
$$
By definition $u\in \bbf(\ob)$ if $u$ is $[-\infty, \infty)$-valued and upper semi-continuous 
on $\ob$ and $u\bigr|_{\O}\in \bbf(\O)$.

The proof of our main result here follows  (as closely as possible) the existence proof 
for the Dirichlet Problem given in [\DDR].

To begin we consider the following analogues of the above.  Let
$$
 E(\ob) \ \equiv\ \left \{ u\in C(\ob): u\bigr|_{\O} \in E(\O)\right\}
\eqno{(\HH.4)}
$$
denote the space of edge functions on $\ob$, and consider the family of edge functions
$$
\cf_{E}(\vf) \ \equiv\ \left \{h\in  E(\ob) : h\bigr|_{\bo}\leq\vf\right\}.
\eqno{(\HH.5)}
$$
A natural question to ask is:

\Qu 1  When is the envelope 
$$
U_E(x) \ \equiv\ \sup_{h\in \cf_{E}(\vf)} h(x) \ \ \text{equal to the solution $H$ defined by (\HH.2)?}
$$

There are two interesting extreme cases where the answer is positive.

\Ex {\HH.1. ($\bbf \equiv \cp$)} Here $E(\O) \equiv {\rm Aff}(\rn)$, the space of affine functions on $\rn$.
In this case 
$$
U_{\rm Aff} \ =\ H_\cp
$$
because, by the Hahn-Banach Theorem, for each point $x_0\in \O$, there exists an affine function $h$
with $h\leq H_\cp$ on $\ob$ and $h(x_0) = H_\cp((x_0)$.

\Ex {\HH.2. ($\bbf \equiv \D$)} Here $E(\O) \equiv \{h\in C(\ob) : h  \bigr|_{\O} \ \text{is $\D$-harmonic}\}$.
Therefore, $H\in \cf_E(\vf)$, proving that
$$
U_\D \ =\ H_{\D}.
$$

For other cases Question 1 remains open, so it is appropriate to consider larger families than $\cf_E(\vf)$.
First, set
\def\emax{{E^{\rm max}}}
\def\elocmax{{E^{\text{loc-max}}}}
$$
\emax(\O) \ \equiv \ \{M : M= \max\{h_1, ... , h_N\}\ \ {\rm with} \ \ h_1, ... , h_N \in E(\O)\}
\eqno{(\HH.6)}
$$
and consider the family
$$
\cf_\emax(\vf) \ \equiv\ \{ M \in \emax(\ob) : M \rest {\bo}\leq\vf\}
\eqno{(\HH.7)}
$$
where by definition $M\in \emax(\ob)$ if $M\in \USC(\ob)$ and $M\rest\O \in \emax(\O)$.
Since the conditions $M \equiv \max\{h_1, ... , h_N\} \in \emax(\ob)$ and $M\rest{\bo}\leq\vf$ 
imply that each $h_k\in \cf_E(\vf)$, we have
$$
U_{\emax} \ =\ \sup_{M\in \cf_\emax(\vf)} M \ =\  \sup_{h\in \cf_E(\vf)} h \ =\ U_E.
\eqno{(\HH.8)}
$$
In particular, $U_{\emax} = H \ \iff\ U_E=H$.

Now we consider a localized version
$$
\cf_{\elocmax}(\vf) \ =\ \{u\in \elocmax(\ob) : u\rest{\bo}\leq \vf\}
\eqno{(\HH.9)}
$$
where by definition $u\in \elocmax(\ob)$ if $u\in \USC(\ob)$ and for each point $x_0\in\O$,
there exists a neighborhood $B_r(x_0)\ss\O$ such that
$$
u \rest{B_r(x_0)} \ \in\ \emax(B_r(x_0)).
\eqno{(\HH.10)}
$$

\Qu 2  When is the envelope 
$$
U\ \equiv\ U_{\elocmax} \ =\ \sup_{u\in \cf_{\elocmax(\vf)}} u \ \ \text{equal to the solution $H$ in (\HH.2)}?
$$

We can answer this question.

\Theorem{\HH.3}
{\sl
If 
$\bbf=\cp^+$ is a minimal subequation and $\bo$ is smooth and strictly $\bbf$-convex,
then}
$$
U\ =\ H.
$$

\pf Since $\bbf\ss\wt\bbf$ (Thm. \DDD.7),  the strict $\wt\bbf$-convexity of the boundary is automatic.

In what follows we shall shorten $\cf_{\elocmax}(\vf)$ to $\cf(\vf)$.

Note that $\cf(\vf) \ss \cf_\bbf(\vf) \Rightarrow U\leq H \Rightarrow U^*\leq H
\Rightarrow$
$$
U^*\bigr|_{\bo} \ \leq\ \vf, \qquad\text{and we also have}
\eqno{(\HH.11a)}
$$
$$
\vf\ \leq\ U_*\bigr|_{\bo} \ \  \qquad\text{proved at the end.}
\eqno{(\HH.11b)}
$$

\Note{\HH.4}  If $\bo$  is strictly convex, then Example \HH.1 shows that 
$\vf = U_\cp\bigr|_{\bo}$ and $U_{\rm Aff} = U_\cp$.  Since $\cp\ss\bbf$ and ${\rm Aff} \ss  E$,
we have $U_\cp \leq U_E \leq U$.  Hence $U_\cp \leq U_*$, so that (\HH.11b) holds
under strict $\cp$-convexity of $\bo$.

These two properties imply the following.
$$
\text{(Boundary Continuity)} 
\qquad
U_*\tobd \ =\ U\tobd \ =\ U^*\tobd \ =\ \vf.
\qquad\qquad
\eqno{(\HH.11)}
$$

By the ``families bounded above property'' we have
$$
U^* \in \bbf(\ob).
\eqno{(\HH.12)}
$$
Note that $H$ (or $\sup_{\bo} \vf$ if you wish) provides an upper bound for $\cf(\vf)$.

Assume for the moment that:
$$
-U_* \in \wt\bbf(\ob).
\eqno{(\HH.13)}
$$
Then the proof is easily completed as follows.  By (\HH.12) and (\HH.13),
$U^*-U_* \in \wt\cp(\ob)$ is subaffine on $\O$ (see [\DDd]).
Moreover, it is $\geq0$ on $\ob$ and equal to zero on $\bo$.  
Hence, by the (MP) for $\wt\cp$, $U^*-U_*$ vanishes on $\ob$.   That is, 
$$
U_*\ =\ U\ =\ U^* \quad{\rm on}\ \ \ob.
\eqno{(\HH.14)}
$$
This proves that $U$ is $\bbf$-harmonic on $\O$ and equal
to $\vf$ on $\bo$.  By uniqueness for the (DP) this proves that
$U=H$ on $\ob$.  Thus it remains to prove (\HH.11b) and  the following.

\Lemma {\HH.5}  \ \ \ \ $-U_*\bigr|_{\O} \in \wt\bbf(\O)$.

\pf  We follow that bump argument given in the proof of Lemma $\wt F$   
in [\DDR, p.\ 455] as closely as possible.

Suppose $-U_*\bigr|_{\O} \notin \wt\bbf(\O)$.  Then there exists $x_0\in \O$, $\e>0$
and $\psi$, a degree-2 polynomial, satisfying
$$
\begin{aligned}
&(a) \ \ -U_* \ \leq \ \psi -\e|x-x_0|^2\quad{\rm near} \ \ x_0, \ \ {\rm and} \\
&(b) \ \ -U_*(x_0) \ =\ \psi(x_0),  \ \ {\rm and} \\
&(c) \ \ \ \  D^2\psi \notin \wt \bbf.
\end{aligned}
\eqno{(\HH.15)}
$$
Rewrite (a) and (c)  as
$$
\begin{aligned}
&(a)' \ \ -\psi \ \leq \ U_*-\e|x-x_0|^2\quad{\rm near} \ \ x_0, \ \ {\rm and} \\
&(c)'  \ \ \ \  D^2(-\psi) \in \Int \bbf.
\end{aligned}
$$
By the key fact  (\HH.1) above we have that
$$
D^2(-\psi) \ =\ e +P \quad {\rm with} \ \ e\in E \ \ {\rm and}\ \ P>0.
\eqno{(\HH.16)}
$$
Therefore
$$
-\psi \ =\ h + \half\bra{P(x-x_0)}{x-x_0}
\eqno{(\HH.17)}
$$
with $h$ a degree 2 polynomial satisfying
$$
 (i) \ \ D^2 h \ =\ e\and (ii)\ \ h(x_0)\ =\ U_*(x_0).
\eqno{(\HH.18)}
$$
The first part is just the statement that
$$
(i)' \ \ \text{$h$ is an edge function on $\rn$}.
\eqno{(\HH.18)(i)'}
$$
Now by (\HH.17) the inequality (a)$'$ says
$$
h+ \half \bra{P(x-x_0)}{x-x_0} \ \leq\ U_*-\e |x-x_0|^2 \quad{\rm on} \ \ B_{r_2}(x_0).
\eqno{(\HH.19)}
$$
Choose $0< r_1<r<r_2$.  Then by (\HH.19)
$$
\begin{aligned}
h+\d \ &<\ U_*\ \ {\rm on}\ \ B_{r_2}(x_0) -  B_{r_1}(x_0)  \\
\end{aligned}
\eqno{(\HH.20)}
$$
where $\d \equiv \inf_{|x-x_0|=r_1} \half \bra{P(x-x_0)}{x-x_0}$ (or $\d=\e r_1^2$ also works).
For each point $y\in \partial B_r(x_0)$ we have $h(y) +\d < U(y)$ by (\HH.20).  Hence, by the definition  
of the $U=U_{\cf^{\elocmax}}$ given in Question 2,
 there exists $u_y \in \cf(\vf)$ with
$$
h(y) +\d \ <\ u_y(y),
\eqno{(\HH.21)}
$$
and since $h$ and $u_y$ are continuous, this holds in a neighborhood of $y$.  Therefore, by compactness,
there exist $u_1, ... , u_N \in\cf(\vf)$ with
$$
h + \d \  < \ u \ \equiv \  \max\{u_1, ... , u_N\} \ \ \text{in a neighborhood of $\partial B_r(x_0)$}.
\eqno{(\HH.22)}
$$
Since $\cf(\vf)$ is closed under taking  the maximum of a finite number of elements, we have
$$
h + \d \  < \ u \   \text{in a neighborhood of $\partial B_r(x_0)$ with $u\in \cf(\vf)$}.
\eqno{(\HH.23)}
$$
This implies that
$$
u' \ \equiv\ 
\begin{cases}
\qquad u \qquad\qquad {\rm on} \ \ \ob - B_r(x_0)  \\
\max\{u, h+\d\}   \quad {\rm on}  \ \ 
   \overline{B_r(x_0)}
\end{cases}
\eqno{(\HH.24)}
$$
is an element of $\cf(\vf)$. (Note that $h+\d$ and hence $u'$ is not necessarily  an element of $\cf_E^{\rm max}(\vf)$.)
Since $u'\in \cf(\vf)$, we have $u'\leq U$ on $\ob$.
In particular, $h+\d \leq U$ on $B_r(x_0)$.  Since $h$ is continuous, this implies $h+\d \leq U_*$, and hence
$$
h(x_0) + \d \ \leq\ U_*(x_0),
\eqno{(\HH.25)}
$$
which contradicts (\HH.18 b)  that $h(x_0)=U_*(x_0)$.

It only remains to do the following.

\noindent
{\bf Proof of (\HH.11b).}
We fix $x_0\in \bo$, and let  $\rho$ be a smooth, strictly $\bbf$-convex 
defining function for $\partial\O$ defined in a neighborhood of $x_0$. 
Then by (\HH.1)  there exist  $\e>0$ and $r>0$ such that
$$
D^2_x \rho -\e I \  \in \  \Int\bbf \ =\ E+\Int \cp  \qquad \forall \, x\in B_r(x_0).
$$
In particular,  
$$
D^2_{x_0} \rho -\e I \ =\  A+ P  \qquad {\rm for}\ \ A\in E \ \ {\rm and}\ \  P>0.
$$
  By adding a linear function to $\half\bra{Ax}x$, we get a quadratic $\psi$ with $D^2_{x_0}\psi = A$ and $\psi(x_0)=0$ so that
$$
\rho(x) - {\e\over 2} |x-x_0|^2  \ =\ \psi(x) +\half \bra {P(x-x_0)}{x-x_0} + O(|x-x_0|^3).
$$
Taking $\e$ smaller, we can get a smaller $r>0$ so that
$$
\rho(x) - {\e\over 2} |x-x_0|^2  \ >\ \psi(x) +\half \bra {P(x-x_0)}{x-x_0} \qquad{\rm for}\ \ x\in \overline{B_r(x_0)} -\{x_0\}.
$$

Since $\rho\leq0$ on $\ob$ we have
$$
- {\e\over 2} |x-x_0|^2 -\half \bra {P(x-x_0)}{x-x_0}  \ \geq\ \psi(x) 
\qquad{\rm for}  \ \ x\in  \overline{B_r(x_0)}\cap\ob.  
\eqno{(\HH.26)}
$$

We now fix $\d>0$ and shrink $r>0$ so that 
$$
\vf(x_0)-\d \ < \ \vf 
\qquad{\rm for}  \ \ x\in  \overline{B_r(x_0)}\cap\bo. 
\eqno{(\HH.27)}
$$

From (\HH.26) above we have that there exists $\eta$ with
$$
0 \ >\ \eta \ \geq \ \psi(x) 
\qquad{\rm for}  \ \ x\in \left( \overline{B_r(x_0)} - B_{r/2}(x_0)\right) \cap\overline\O. 
\eqno{(\HH.28)}
$$

We now  consider the edge function
$$
\Psi(x) \ \equiv\ \vf (x_0)-\d + C\psi(x).
\eqno{(\HH.29)}
$$
By (\HH.28) we see that for  $C > > 0$ we will have
 $$
 \Psi(x) \ < \  \inf\vf \qquad {\rm on}\ \ \left( \overline{B_r(x_0)} -  B_{r/2}(x_0)\right)\cap\O
 $$
Therefore
$$
\underline u \ \equiv \ 
\begin{cases}    
\inf_{\bo} \vf \qquad {\rm on} \ \ \overline{ \O} - B_{r/2}(x_0)  \\
\max\{\Psi, \inf_{\bo}\vf\} \ \ {\rm on} \ \ \overline{B_r(x_0)}\cap\ob
\end{cases}
$$
is a well defined function on $\ob$, and it  is locally   the maximum
of edge functions.  Furthermore, by (\HH.27) and (\HH.29) we see that $\underline u \leq\vf$ on $\bo$.
Hence,  $\underline u$ is in our Perron family for the Dirichlet problem,
and so we have 
$
\underline u  \leq U,
$
which  implies that 
$
\underline u  \leq U_*.
$
In particular, $\underline u (x_0) = \vf(x_0)-\d \leq U_*(x_0)$.
Taking $\d\to 0$ shows that $\vf(x_0)\leq U_*(x_0)$.
This proves (\HH.11b) and therefore Theorem \HH.3.
\qed

\vskip.4in


\def\item{}

\centerline
{
\headfont  References
}

\noindent
\item{[\Bry]}  R. Bryant and   F. R. Harvey, {\sl Submanifolds in Hyper-Kahler Geometry},
J. A. M. S., {\bf 2} no. 1 (1989),  1-31.

\smallskip

\noindent
\item{[\CRA]}   M. G. Crandall,  {\sl  Viscosity solutions: a primer},  
pp. 1-43 in ``Viscosity Solutions and Applications''  Ed.'s Dolcetta and Lions, 
SLNM {\bf 1660}, Springer Press, New York, 1997.

 \smallskip

\noindent
\item{[\CIL]}   M. G. Crandall, H. Ishii and P. L. Lions {\sl
User's guide to viscosity solutions of second order partial differential equations},  
Bull. Amer. Math. Soc. (N. S.) {\bf 27} (1992), 1-67.

\smallskip

 \noindent
\item{[\Har]}   F. R. Harvey, Spinors and Calibrations, Perspectives in Mathematics, 9. Academic Press, Inc., Boston, MA, 1990.

\smallskip

 \noindent
\item{[\DDd]}   F. R. Harvey and H. B. Lawson, Jr., {\sl  Dirichlet duality and the non-linear Dirichlet problem},    Comm. on Pure and Applied Math. {\bf 62} (2009), 396-443. ArXiv:math.0710.3991.

\smallskip

 \noindent
\item{[\PUP]}  \ \----------,  {\sl  Plurisubharmonicity in a general geometric context},  
Geometry and Analysis {\bf 1} (2010), 363-401. ArXiv:0804.1316

 \noindent
\item{[\DDR]}  \ \----------, {\sl Dirichlet Duality and the Nonlinear Dirichlet Problem on Riemannian Manifolds},  J. Diff. Geom. {\bf 88} (2011), 395-482.   ArXiv:0912.5220.

  \smallskip

 \noindent
\item {[\Survey]}  \ \----------,   {\sl  Existence, uniqueness and removable singularities
for nonlinear partial differential equations in geometry},\ 
 pp. 102-156 in ``Surveys in Differential Geometry 2013'', vol. 18,  
H.-D. Cao and S.-T. Yau eds., International Press, Somerville, MA, 2013.
ArXiv:1303.1117.

\smallskip

\noindent
\item  {[\AE]}   \ \----------,    {\sl The AE Theorem and Addition Theorems for quasi-convex functions,} 
 \ ArXiv:1309:1770.

\smallskip

\noindent
\item  {[\LA]}   \ \----------,     {\sl  Lagrangian potential theory    and a Lagrangian equation of Monge-Amp\`ere type.}

\smallskip

\end{document}